\documentclass[english]{smfart}

\usepackage{amsmath,amssymb}
\usepackage{enumerate}

\numberwithin{equation}{section}
\newtheorem{lemma}[equation]{Lemma}
\newtheorem{problem}{Problem}
\newtheorem{question}{Question}
\newtheorem{proposition}[equation]{Proposition}
\newtheorem{theorem}{Theorem}
\newtheorem{corollary}[equation]{Corollary}

\newtheorem{fact}[equation]{Fact}
\newtheorem{definition}[equation]{Definition}
\newtheorem{remark}[equation]{Remark}

\newtheorem{exercise}[equation]{Exercise}

\newenvironment{demo}{\medskip\noindent{\sc
Proof:}}{\hfill$\square$\medskip}

\newenvironment{demof}[1]{\medskip\noindent{\sc
Proof of #1:}}{\hfill$\square$\medskip}

\newcommand{\new}[1]{{\bf #1}}
\newcommand{\ignore}[1]{}
\newcommand\step[2]{\medbreak\noindent{{\bf Step #1:} {\it
#2}\medbreak}}

\newcommand\diam{{\operatorname{diam}}}
\newcommand\eps{\epsilon}
\newcommand\fol{{\operatorname{fol}}}
\newcommand\HS{{\hat\Sigma}}

\newcommand\INT{{\operatorname{int}}}

\newcommand\Lip{{\operatorname{Lip}}}
\newcommand\loc{{\operatorname{loc}}}

\newcommand\MD{{\mathcal D}}
\newcommand\NN{{\mathbb N}}

\newcommand\PMM{{\operatorname{PMM}}}
\newcommand\Prob{{\operatorname{Prob}}}
\newcommand\RR{{\mathbb R}}

\newcommand\supp{{\operatorname{supp}}}

\renewcommand\top{{\operatorname{top}}}

\newcommand\ZZ{{\mathbb Z}}

\begin{document}

\title{A Minicourse on Entropy Theory on the Interval}

\begin{abstract}
We give a survey of the entropy theory of interval maps as it can
be analyzed using ergodic theory, especially measures of maximum
entropy and periodic points. The main tools are (i) a suitable
version of Hofbauer's Markov diagram, (ii) the shadowing property
and the implied entropy bound and weak rank one property, (iii)
strongly positively recurrent countable state Markov shifts.
Proofs are given only for selected results. This article is based
on the lectures given at the {\sl Ecole th\'ematique de th\'eorie
ergodique} at the C.I.R.M., Marseilles, in April 2006.
\end{abstract}

\author{J\'er\^ome Buzzi}

\address{Centre de Math\'ematiques\\
Ecole polytechnique\\91128 Palaiseau Cedex\\France}

\email{buzzi@math.polytechnique.fr}

\urladdr{www.jeromebuzzi.com}

\thanks{J.B. wishes to thank the organizers of the {\it Ecole th\'ematique de th\'eorie ergodique}
at the CIRM in Marseilles.}


\maketitle

\tableofcontents

\newpage

\section{Introduction}

We are going to give a very selective survey of interval dynamics,
mainly (but not exclusively) those defined by maps with finitely
many critical points or discontinuities (see Definition
\ref{def-pmm} below). We focus on ``complexity'' as defined
through entropy as seen from an \emph{ergodic theory} point of
view. Ergodic theory will be for us both a powerful tool and a
guide to the ``right'' questions. In particular, we shall
concentrate on aspects not dealt previously in book form (like the
classic treatise \cite{MS} or \cite{AM} for another, non-ergodic,
point of view on the same subject) around measures of maximum or
large entropy.

As the lectures given in Luminy, these notes are intended to be
accessible to readers with only a basic knowledge of dynamical
systems. From a technical point of view we shall only assume (1)
measure theory (e.g.,  the very first chapters of (\cite{Rudin});
(2) the Birkhoff ergodic theorem, explained in any textbook on
ergodic theory (see, e.g., \cite{Walters}).

We shall deal mainly with the following rather well-behaved class
of one-dimensional dynamics:

\begin{definition}\label{def-pmm}
$I$ will always denote a compact interval of $\mathbb R$. A
self-map of $I$ is
\new{piecewise monotone} if there exists a partition of $I$ into
\emph{finitely many} subintervals (the "pieces") on each of which
the restriction of $f$ is continuous and \emph{strictly} monotone.
Note that one can always subdivide the pieces. The \emph{natural
partition} is the set of interiors (relatively to $\mathbb R$) of
the pieces in such a partition with minimum cardinality.

We denote the set of piecewise monotone maps by $\PMM(I)$ and its
topology is defined by $d(f_1,f_2)<\eps$ if $f_1$ and $f_2$ both
admit natural partitions with $n$ pieces with endpoints
$a^i_j,b^i_j$ such that $|a^1_j-a^2_j|<\eps$, $|b^1_j-b^2_j|<\eps$
and $|f^1(a^1_j+(b^1_j-a^1_j)t)- f^2(a^2_j+(b^2_j-a^2_j)t)|<\eps$
for all $t\in [0,1]$.
\end{definition}

\begin{remark}
A piecewise monotone map is \emph{not} assumed to be continous.
The above topology induce a topology on $C^0(I)$ which is neither
stronger nor weaker than the usual one.
\end{remark}

Let us outline these notes. We shall recall in Section
\ref{sec-entropy-gal} some general facts about entropy for smooth,
topological or probabilistic dynamical systems. In Section
\ref{sec-h-interval}, after recalling the basics of the symbolic
dynamics of piecewise monotone maps, we give combinatorial and
geometric formulations of the entropy on the interval using in
particular the special form of its symbolic dynamics defined by
the kneading invariants. We discuss the continuity and
monotonicity of the entropy function over maps and over invariant
probability measures in Section \ref{sec-function-entropy}.

To get to the global structure we shall use Hofbauer's Markov
diagram explained in Section \ref{sec-diagram}. This will leave a
part of the dynamics which we analyze using "shadowing" in Section
\ref{sec-shadow}. The main part of the dynamics is reduced to a
Markov shift with countably many states but most of the properties
of the finite case (section \ref{sec-SPR}). We apply these tools
to the piecewise monotone maps getting a precise description of
their measures of large or maximum entropy and their peirodic
points, including a complete classification from this point of
view (section \ref{sec-application}). We conclude in Section
\ref{sec-conclusion} by mentioning further works which either
analyze more precisely the piecewise monotone maps or apply the
techniques presented here to more general settings, following the
idea that one-dimensional dynamics should be the gateway to (more)
general results.

\begin{remark}
Theorems, Problems and Questions are numbered consecutively and
independently throughout the paper. All other items are numbered
in a common sequence within each section.

{\bf Exercises} should be rather straightforward and quick
applications of techniques and ideas exposed in the text. {\bf
Problems} are more ambitious projects that (I believe) can be
solved by standard techniques - but have not been done yet, to the
best of my knowledge. {\bf Questions} are problems that I don't
know how to handle.
\end{remark}

\section{Generalities on Entropies}\label{sec-entropy-gal}

The dynamical entropies have a long story and are related to very
basic notions in statistical physics, information theory and
probability theory (see, e.g., \cite{KellerEntropy}).

Putting the work of Shannon on a completely new level of
abstraction, Kolmogorov and Sinai defined in 1958 the
\new{measured entropy}\footnote{We prefer this nonstandard terminology
to the usual, but cumbersome \emph{measure-theoretic} and even
more to the confusing \emph{metric entropy}.} of any endomorphism
$T$ of a probability space $(X,\mu)$ as follows. For any finite
measurable partition $P$, we get a $P$-valued process whose law is
the image of $\mu$ by the map $x\mapsto (P(x),P(Tx),\dots)\in
P^\NN$. This process has a mean Shannon entropy:
 $$
   h(T,\mu,P) := \lim_{n\to\infty} \frac1n H(\mu,P^n)
 $$
where $H(\mu,Q)=\sum_{A\in Q} -\mu(A)\log\mu(A)$ and
 $$
   P^n := \{ \left<A_0\dots A_{n-1}\right>:=A_0\cap T^{-1}A_1\cap\dots T^{-n+1}A_{n-1}\ne\emptyset
    : A_i\in P\}
 $$
The elements of $P^n$ are called the (geometric)
\new{$P,n$-cylinders}.

The Kolmogorov-Sinai $h(T,\mu)$ is then the supremum of the
Shannon entropies of all processes over finite alphabet
"contained" in the considered dynamical system. We refer to the
many excellent texts (see, e.g.,
\cite{Mane,Sinai1,Sinai2,Petersen,Rudolph, Shields}) for more
information and only quote a few facts here.

\medbreak

This supremum can look forbidding. However, Sinai showed: if
$P_1,P_2,\dots$ is an increasing\footnote{That is, the elements of
$P_n$ are union of elements of $P_{n+1}$ for all $n$.} sequence of
finite measurable partitions such that $\{T^{-k}P_n:k,n\in\NN\}$
generates the $\sigma$-algebra of measurable subsets of $X$, then:
 $$
    h(T,\mu) = \sup_{n\geq1} h(T,\mu,P_n)
 $$
Note the case where all $P_n$'s are the same partition (said then
to be \new{generating} under $T$).

The measured entropy is an invariant of measure-preserving
conjugacy: if $(X,T,\mu)$ and $(Y,S,\nu)$ are two measure
preserving maps of probability spaces and if $\psi:X\to Y$ is a
bimeasurable bijection $\psi:X\to Y$ of probability spaces such
that $\nu=\mu\circ\psi$ and $\psi\circ T=S\circ\psi$ (i.e., $\psi$
is an isomorphism of $(X,T,\mu)$ and $(Y,S,\nu)$) then $T$ and $S$
have the same entropy.

Stunningly, this invariant is complete for this notion of
isomorphism among Bernoulli automorphisms   according to
Ornstein's theory (see \cite{Petersen} for an introduction and
\cite{Rudolph} for a complete treatment). Ornstein theory also
(and perhaps more importantly) shows that many natural systems are
measure-preserving conjugate to such a system (see \cite{OW}).

Let us note that the above define the entropy wrt not necessarily
ergodic invariant probability measure. One shows also that
$h(T,\mu)$ is an \emph{affine function} of $\mu$ so that, if
$\mu=\int \mu_x \nu(dx)$ is the ergodic decomposition of $\mu$,
then $h(T,\mu)=\int h(T,\mu_x) \, \nu(dx)$. We also note that
$h(T^n,\mu)=|n|h(T,\mu)$ for all $n\geq1$ (for all $n\in\ZZ$ if
$T$ is invertible).

From our point of view, the real meaning of measured entropy is
given by the Shannon-McMillan-Breiman theoreom:

\begin{theorem}[Shannon-McMillan-Breiman]\label{thm-SMB}
Let $T$ be a map preserving a probability measure $\mu$ on a space
$X$. Assume that it is ergodic. Let $P$ be a finite measurable
partition of $X$. Denote by $P^n(x)$ the set of points $y$ such
that $f^kx$ and $f^ky$ lie in the same element of $P$ for $0\leq
k<n$. Then, as $n\to\infty$:
 $$
    \frac1n \log \mu( P^n(x) ) \to h(T,P,\mu) \qquad \text{a.e. and
    in }L^1(\mu)
 $$
\end{theorem}

\begin{corollary}\label{coro-SMB}
Let $T$ and $P$ be a as above. For $0<\lambda<1$, let
$r(P,n,\mu,\lambda)$ be the minimum cardinality of a collection of
$P,n$-cylinders the union of which has $\mu$-measure at least
$\lambda$. Then:
 $$
    h(T,\mu,P) = \lim_{n\to\infty} \frac1n\log r(P,n,\mu,\lambda)
 $$
\end{corollary}

\begin{exercise}
Consider $X=\{0,1\}^\NN$ together with the shift $\sigma$ and the
product probability $\mu$ induced by $(p,1-p)$. Show that:
 $$
   h(\sigma,\mu) = -p\log p-(1-p)\log(1-p)
 $$
\end{exercise}

The following general fact is especially useful in dimension $1$:

\begin{theorem}[Rokhlin formula]\label{theo-rokhlin}
Let $T$ be an endomorphism of a probability space $(X,\mu)$.
Assume that there is a generating countable measurable partition
of $X$ into pieces $X_i$ on each of which the restriction $T|X_i$
is one-to-one and non-singular wrt $\mu$ and
$-\sum_{i}\mu(X_i)\log\mu(X_i)<\infty$. Let $JT=d\mu\circ
(T|X_i)/d\mu$ if $x\in X_i$. Then:
 $$
   h(T,\mu) = \int \log JF \, d\mu
 $$
\end{theorem}

Compare this with Pesin's formula (below).

\begin{exercise}
Recover the result of the previous exercise by applying Rokhlin
formula to a suitable dynamical system on the interval.
\end{exercise}

\subsection{Bowen-Dinaburg formula}

R. Adler, A. Konheim,  and M. McAndrew defined in 1965 a
topological counterpart of the measured entropy for arbitrary
continuous maps of compact spaces. Where measured entropy counts
orbit segments "typical" for the measure (Corollary
\ref{coro-SMB}), topological entropy counts all orbit segments.

The original definition of topological entropy was given in terms
of open covers. In the case of a compact metric space $(X,d)$, we
shall use the Bowen-Dinaburg formula which we now explain.

\begin{definition}
For a real number $\eps>0$ and a positive integer $n$, the
$(\eps,n)$-ball at some $x\in X$ is:
 $$
   B(\eps,n,x) := \{y\in X:\forall 0\leq k<n\; d(T^ky,T^ky)< \eps
   \}
 $$
The $(\eps,n)$-covering number of $Y\subset X$ is:
 $$
   r(\eps,n,Y) := \min \{\# C: Y\subset \bigcup_{y\in C}
   B(\eps,n,y) \}
 $$
The entropy of $Y$ is then:
 $$
   h(T,Y):= \lim_{\eps\to 0^+} h(T,Y,\eps) \text{ with }
     h(T,Y,\eps) := \limsup_{n\to\infty} \frac1n \log r(\eps,n,Y)
 $$
The \new{topological entropy} of $T$ is:
 $$
   h_\top(T) := h(T,X).
 $$
\end{definition}

Thus, $h_\top(T)$ counts \emph{all} orbit segments.

\begin{remark}

1. $h_\top(T)$ is an invariant of topological conjugacy.

2. One can replace $\limsup$ by $\liminf$ in the definition of
$h(T,Y,\eps)$. The number $h_\top(T)=h(T,X)$ is then left
unchanged.
\end{remark}

\begin{question}
It is unknown if for $X=Y$ the $\limsup$ in $h(T,X,\eps)$ is in
fact a limit.
\end{question}

\begin{exercise}
Show that if $T$ is Lipschitz with constant $L$ on a compact
manifold of dimension $d$ then $h_\top(T)\leq d\log^+ L$. {\it
Note:} A very important refinement is Ruelle-Margulis inequality
(see \ref{thm-ruelle-inequality} and its descendants. It is false
if $T$ is only piecewise continuous and Lipschitz, see
\cite{Affine,Kruglikov}.
\end{exercise}

\begin{exercise}
Let $\sigma:\mathcal A^\NN\to\mathcal A^\NN$ be the shift map:
$(\sigma(A))_n=A_{n+1}$. Let $\Sigma$ be a \new{subshift} on a
finite alphabet $\mathcal A$, i.e., $\Sigma$ is a closed,
$\sigma$-invariant subset of $\mathcal A^\NN$. Show that
 $$
    h_\top(\sigma|\Sigma) = \lim_{n\to\infty}\frac1n\log
      \#\{w\in \mathcal A^n:\exists x\in\Sigma\text{ s.t. }
      x_0\dots x_{n-1}=w\}
 $$
If $\Sigma$ is the \new{subshift of finite type} (see \cite{LM}
for background) defined by a matrix $A:\mathcal A\times\mathcal
A\to\{0,1\}$ according to $x\in\Sigma$ iff $x\in\mathcal A^\NN$
and, for all $n\in\NN$, $A(x_n,x_{n+1})=1$, then
 $$
   h_\top(\sigma|\Sigma)=\log \rho(A)
 $$
where $\rho(A)$ is the spectral radius of $A$.
\end{exercise}

\subsection{Katok's formula}

Let $T$ be an endomorphism of a probability space $(X,\mu)$.
Assume that $X$ is also a compact metric space (the measurable
structure being the Borel one). Katok \cite{Katok0} observed that
one can compute the measured entropy by a variant of the
Bowen-Dinaburg formula:

\begin{definition}
The $(\eps,n)$-covering number of a probability $\mu$ wrt a
parameter $\lambda\in(0,1)$ is:
 $$
   r(\eps,n,\mu,\lambda) := \min \left\{\#: \mu\left(\bigcup_{y\in
   C}B(\eps,n,y)\right) >\lambda \right\}
 $$
\end{definition}

\begin{proposition}[Katok]
Assume that $(T,\mu)$ is ergodic. Then the measured entropy of
$(T,\mu)$ is equal to:
 $$
   h(T,\mu) = \lim_{\eps\to0^+} h(T,\mu,\eps) \text{ with }
   h(T,\mu,\eps) = \limsup_{n\to\infty} \frac1n\log r(\eps,n,\mu,\lambda)
 $$
for any $0<\lambda<1$. $h(T,\mu,\eps)$ is independent of
$\lambda$.
\end{proposition}

\begin{remark}
One can replace $\limsup$ by $\liminf$ in the above definition.
However it is not known if this is in fact a limit
\end{remark}

\subsection{Ruelle-Margulis Inequality}

Let $f:M\to M$ be a self-map of a compact manifold with an ergodic and
invariant probability measure $\mu$ and an open subset $U$ such that $\mu(U)=1$
and $f|U$ is $C^1$.

The Lyapunov exponents of $(f,\mu)$ are numbers $\lambda_1\geq\lambda_2\geq
\dots\lambda_d\geq-\infty$ such that, for a.e. $x\in M$:
  $$
   \dim\overline{\{v\in T_xM:\lambda(x,v)=\lambda\}} = \#\{i:\lambda_i=\lambda\}.
  $$
Oseledets Theorem ensures that such a collection of $d$ numbers
exists if $\log^+\|f'\|$ is $\mu$-integrable (see, e.g.,
\cite{Krengel}).

\begin{theorem}[Ruelle-Margulis inequality]\label{thm-ruelle-inequality}
If $\|f'\|$ is bounded over $U$, then
 $$
    h(f,\mu) \leq \sum_{i=1}^d \lambda_i^+
  $$
\end{theorem}

See \cite{KatokStrelcyn} for refinements (i.e., how fast can one
let $f$ blow up near $\partial U$).

\subsection{Variational Principle}

The above formulas imply immediately that $h(T,\mu)\leq h_\top(T)$
for all measurable maps $T$ of a compact metric space and all
invariant probability measures $\mu$. In fact much more is true:

\begin{theorem}[Goodman, Dinaburg]
For a continuous map of a metric space, the variational principle
holds:
 $$
   h_\top(T) = \sup_\mu h(T,\mu)
 $$
The supremum can be taken either over all invariant probability
measures; or over all ergodic and invariant probability measures.
\end{theorem}

The classical proof is due to Misiurewicz (see the textbooks,
e.g., \cite{Petersen})

\begin{remark}

1. This of course "justifies" the definition of topological
entropy to the extent that this quantity is shown to depend only
on the measurable structure of $T:X\to X$ and not on its
topological (or metric) structure.

2. It was not at all obvious that one could find measures with
"almost all the topological complexity" of the system. For
instance, the opposite is true in the setting of Birkhoff theorem:
e.g., for the shift on $\{0,1\}^\NN$, the set of points for which
$\lim_{n\to\infty}\frac1n\sum_{k=0}^{n-1} \phi\circ \sigma^k(x)$
exists for all continuous function $\phi$ is of full measure but
is of first Baire category (i.e., it is negligible from the
topological point of view).
\end{remark}

\begin{exercise}
Prove this last assertion.
\end{exercise}

The variational principle leads to the study of invariant
probability measures that achieve the topological entropy. More
generally, one makes the following

\begin{definition}\label{def-max-meas}
For a Borel map $T:X\to X$, the entropy is defined as
 $$
   h(T) := \sup_\mu h(T,\mu)
 $$
where $\mu$ ranges over all the invariant probability measures.
The \new{maximum measures} are the ergodic, invariant probability
measure maximizing entropy, i.e., such that $h(T,\mu)=h(T)$.
\end{definition}

Recall that the entropy is an affine function of the measure,
hence the invariant probability measures maximizing entropy are
exactly the closed convex hull of the ergodic ones, i.e., the
maximum measures.

\medbreak

Maximum measures sometimes do not exist and sometimes are not
unique. For instance, for each finite $r$, there are $C^r$ maps
with no maximum measures and also $C^r$ maps with infinitely many
maximum measures \cite{BuzziSIM,RuetteMax} (see
\cite{MisiurewiczNo} for examples of diffeomorphisms). But we
shall see that for many nice systems, such measures exist, have
some uniqueness and describe important characteristics of the
dynamics, e.g., the repartition of the periodic points.

\subsection{Misiurewicz's local entropy}

A key aspect of measured entropy is the uniformity (or lack
thereof) in the limits, when $\mu$ ranges over the invariant
probability measures, $h(f,\mu)=\lim_{\eps\to0} h(f,\mu,\eps)$. It
corresponds to the existence of "complexity at arbitrarily small
scales". The simplest way\footnote{Much deeper results have been
obtained in the study of symbolic extension entropy, see the
references in Section \ref{sec-conclusion}.} to measure this is
the following quantity introduced by Misiurewicz (under the name
\emph{topological conditional entropy}) \cite{Misiu-hloc}:

\begin{definition}[Misiurewicz]
The \new{local entropy} of a self-map $f$ of a metric space
$(X,d)$ is:
 $$
   h_\loc(f) := \lim_{\eps\to 0} h_\loc(f,\eps)
   \text{ with }
   h_\loc(f,\eps) := \lim_{\delta\to0} \limsup_{n\to\infty} \sup_{x\in X}
     \frac1n\log r(\delta,n,B(\eps,n,x)).
 $$
\end{definition}

\begin{exercise}
Show that, if $X$ is compact and $f$ is continuous, then
$h_\loc(f)=\lim_{\eps\to 0} \tilde h_\loc(f,\eps)$ where $\tilde
h_\loc(f,\eps) := \sup_{x\in X} \lim_{\delta\to0}
\limsup_{n\to\infty} \frac1n\log r(\delta,n,B(\eps,n,x))$.
\end{exercise}

\begin{exercise}\label{exo-hloc-iter}
Show that if $X$ is compact, but $f$ not necessarily continuous,
then $h_\loc(f^n)=|n| h_\loc(f)$ for all $n\in\NN$ (or $n\in\ZZ$
if $f$ is invertible).
\end{exercise}

A simple but important consequence motivating the above definition
is:

\begin{proposition}\label{prop-hloc}
$h_\loc(f)$ bounds the defect in upper semicontinuity of
$\mu\mapsto h(f,\mu)$ over the set of invariant probability
measure with the weak star topology: for any sequence
$\mu_n\to\mu$,
 $$
   \limsup_{n\to\infty} h(f,\mu_n) \leq h(f,\mu)+h_\loc(f)
 $$
In particular, if $h_\loc(f)=0$, then the above map is upper
semi-continuous and therefore achieves its maximum: there exists
at least one maximum entropy measure.
\end{proposition}

\begin{exercise}
Prove this result. \emph{Hint:} for any $0<\eps'<\eps$, for all
$n$ large enough, $h(f,\mu_n,\eps)\leq h(f,\mu,\eps')$.
\end{exercise}

$h_\loc(f)=0$ for $C^\infty$ maps follows from Yomdin's theory to
which we now turn.

\subsection{Yomdin's theory}

Yomdin's theory analyzes the local complexity of differentiable
maps. It explains why $C^\infty$ interval maps are so much like
piecewise monotone maps (see \cite{BuzziSIM} or section
\ref{sec-conclusion}).

\begin{theorem}[following Yomdin]\label{theo-yomdin}
Let $f:\RR^d\to\RR^d$ be a $C^r$ map with all partial derivatives
up to order $r$ bounded by $1$, except for the derivatives of
order $1$ maybe.

Let $\sigma:[0,1]^k\to\RR^d$ extend to a $C^r$ map defined on a
neighborhood of $[0,1]^k$ with all partial derivatives up to order
$r$ bounded by $1$.

Then there exists a tree $\mathcal T$ (i.e., a collection of
finite sequences such that if $\tau_1\dots\tau_n\in\mathcal T$
then $\tau_1\dots\tau_k\in\mathcal T$ for all $1\leq k\leq n$) and
a family of maps $\phi_\tau:[0,1]^k\to[0,1]^k$, $\tau\in\mathcal
T$, with the following properties. Set for any
$\tau=\tau_1\dots\tau_n\in\mathcal T_n$ (the subset of sequences
of length $n$),
 $$\Phi_{\tau_1\dots\tau_n}:=\phi_{\tau_1\dots
\tau_n}\circ\phi_{\tau_1\dots
  \tau_{n-1}}\circ \dots\circ\phi_{\tau_1}.
  $$
Then, for all  $n\geq1$,
 \begin{itemize}
  \item $\sigma^{-1}\left(B(1,n,\sigma(0))\right)
   \subset \bigcup_{\tau\in\mathcal T_n} \Phi_\tau([0,1]^k)
   \subset \sigma^{-1}\left(B(2,n,\sigma(0))\right)$
  \item for each $\tau\in\mathcal T_n$, $\phi_\tau$, $\Phi_\tau$ and
  $f^{n}\circ\sigma\circ\Phi_\tau$
  have all their partial derivatives up to
  order $r$ bounded by $1$;
  \item $\#\mathcal T_n\leq C(r,d)\Lip^+(f)^{n(k/r)+1}$ where $\Lip^+(f)$ is the
  maximum of the Lipschitz constant of $f$ and $1$.
  \end{itemize}
\end{theorem}

A slightly weaker claim was prove by Yomdin \cite{Yomdin}. However
simple modifications (made in \cite{BuzziSIM}) yields the above
statement. We do not give a proof here, but only sketch the main
ideas.

\subsubsection*{Sketch of Proof}

The first step of the proof is to observe that one can approximate
$f^n\circ\sigma:[0,1]^k\to\mathbb R^d$ by its Taylor polynomial,
after restricting it to a cube of linear size
$o(\Lip^+(f)^{-n/r})$. There are essentially $\Lip^+(f)^{n(k/r)}$
such cubes in the domain $[0,1]^k$. This accounts for the
corresponding factor in the above Theorem.

The second step of the proof considers the intersection with the
unit cube of the graph of the polynomial approximation of
$f^n\circ\sigma$. This is a \new{semi-algebraic subset}, i.e., a
subset of $\mathbb R^d$ defined by finitely many equalities and
inequalities involving only polynomials. It turns out that such
subsets can be written as the union of the images of a
\emph{constant number} of $C^r$ maps with all partial derivatives
bounded by $1$. Here by a constant number, we mean one that
depends only on the order of differentiability $r$, dimension $d$
and the list of the degrees of the polynomials defining the
subset. This natural but deep fact was explained in \cite{Gromov}
and a complete proof can be found in \cite{Burguet}.

\subsubsection*{Applications}

The initial motivation of Yomdin was the following theorem:

\begin{theorem}[Yomdin]
Let $f$ be a $C^r$ self-map of a compact $d$-dimensional
Riemannian manifold $M$. Define the volume growth of $f\in
C^r(M)$, $r\geq1$, by:
 $$
   v(f) := \max_{0\leq k\leq d}\sup_S \limsup_{n\to\infty}
      \frac1n \log \int_S \|\Lambda^k(f^n)'(x)\| \, dvol_S
 $$
where $S$ ranges over all compact subsets of $k$-dimensional $C^r$
submanifolds of $M$, $vol_S$ is the volume induced by the
restriction to $S$ of the Riemannian structure of $M$ and
$\|\Lambda^k (f^n)'(x)\|$ is the Jacobian corresponding to $vol_S$
and $vol_{f^nS}$. Thus $v(f)$ is the growth rate of all volumes
with multiplicity.

Then,
 \begin{equation}\label{eq-vol-leq-h}
    v(f) \leq h_\top(f) + \frac{d}{r}\log^+\Lip^+(f)
 \end{equation}
In particular, for $C^\infty$ maps, the logarithm of the spectral
radius of the action of $f$ on the total homology is bounded by
the topological entropy.
\end{theorem}

\begin{exercise}
Deduce the above theorem from Theorem \ref{theo-yomdin}.
\emph{Hint:} For any $0\leq k\leq n$, bound $\int_B
\|\Lambda^k(f^n)'(x)\| \, dvol_B$ for $B$ contained in an
$(\eps,n)$-ball.
\end{exercise}

{\it Note:} The last assertion is \emph{Shub conjecture} in the
case of $C^\infty$ maps. The conjecture remains open for less
regular functions, in particular $C^1$ functions (for which it was
stated).

The bound (\ref{eq-vol-leq-h}) is sharp among maps, see
\cite{Yomdin}.

\begin{exercise}
Prove that for a $C^r$ map of a compact $d$-dimensional manifold,
$h_\loc(f)\leq \frac{d}{r}\log\Lip^+(f)$. \emph{Hint:} apply
Theorem \ref{theo-yomdin} to a suitable $F$ representing a scaled
version of $f$ and $\sigma$ being a suitable chart and consider
$\bigcup_{w\in\mathcal T_n} \sigma\circ\Phi_\tau(Q_\eps)$ where
$Q_\eps$ is an $\eps$-dense subset of $[0,1]^d$.
\end{exercise}

Even more interestingly (for higher dimensional dynamics),
Yomdin's theory implies a sort of submultiplicative property among
families of maps of the form $f^n\circ\sigma$ where $f$ is the
dynamics and $\sigma$ is a parametrized disk (see
\cite{BuzziSMF}).

\section{Computing the Entropy on the Interval}\label{sec-h-interval}

We now turn to one-dimensional dynamics. Here entropy is somewhat
"explicit". This first appears in two results of Misiurewicz which
say that the entropy for interval maps is essentially a
combinatorial phenomenon.

\subsection{Symbolic Dynamics}

We first recall the general definition of the symbolic dynamics of
a piecewise map before explaining the classical description of the
subshifts thus obtained in the class of piecewise monotone maps.

\begin{definition}
Let $f$ be a self-map of a Baire topological space $X$ with a
collection $P$ of disjoint open subsets with dense union. The
symbolic dynamics of $f$ wrt $P$ is the left-shift $\sigma$ acting
on:
 $$
    \Sigma_+(f,P) := \overline{\{\iota(x):x\in X'\}} \subset P^\NN
 $$
where the topology on $P^\NN$ is the product of the discrete
topology of $P$ and $X':=\bigcap_{k\geq0}f^{-k}\bigcup_{A\in P}A$
and
 \begin{equation}
    \iota(x) := A \in P^\NN \text{ such that } \forall k\geq0\; f^kx\in
    A_k.
 \end{equation}
$\iota(x)$ is called the \new{itinerary} of $x$.
\end{definition}

\begin{exercise}
For $A_0\dots A_{n-1}\in P^n$, let $$\left<A_0\dots
A_n\right>:=A_0\cap f^{-1}A_1\cap\dots \cap f^{-n+1}A_{n-1}\subset
X.$$ Check that $\Sigma_+(f,P) = \{A\in P^\NN: \forall n\geq0\;
\left<A_0\dots A_n\right>\ne\emptyset \}$.
\end{exercise}

For piecewise monotone map, the symbolic dynamics is very close to
the interval dynamics:

\begin{exercise}\label{exo-sd-pmm}
Show that for a piecewise monotone map $f$ with its natural
partition $P$, $\iota$ is well-defined except on a countable set.
Show that $\iota$ is one-to-one except on a set negligible wrt any
nonatomic invariant probability measure. \emph{Hint:} Consider the
homtervals of $f$, i.e., the maximum open intervals on which $f^k$
is monotone and continuous for all $k\geq1$.

Give a counter-example for a smooth interval map. Can you obtain
$h_\top(\sigma|\Sigma_+(f,P))>h_\top(f)$ or the reverse
inequality? See \cite{Affine,Kruglikov} for some piecewise affine
examples in higher dimensions.
\end{exercise}

\subsection{Kneading theory}\label{sec-kneading}

Let $(T,P)$ be a piecewise monotone map.

 Endpoints of pieces in $P$ have no itinerary in the proper
sense. However, one can define left/right itineraries of any point
$x\in I$:
 $$
    \iota(x\pm) := \lim_{\eps\to0^+} \iota(x\pm\eps).
 $$
If $c_0<c_1<\dots<c_N$ are the endpoints of the intervals in the
partition $P$ of $I=[c_0,c_N]$, the set of kneading invariants is:
 $$
    K(T,P):=\{\iota(\gamma):
    \gamma\in\{c_0+,c_1-,c_1+,\dots,c_N^-\}\}
 $$
These are obviously invariants of topological conjugacy. In fact,
they determine completely the symbolic dynamics (hence, they are
(essentially) complete invariants of topological conjugacy -- see
\cite{MS} for the details).

\begin{definition}
$T$ defines a \new{modified lexicographic order} on $P^\NN$
according to:
 \begin{multline*}
    A \preceq_T B \iff A = B \text{ or }
      (A_0\dots A_{n-1}=B_0\dots B_{n-1} \text{ and } \\
        \left<A_0\dots A_{n-1}A_n\right> < \left<B_0\dots B_{n-1}B_n\right>)
 \end{multline*}
where $I<J$ means that for every $x\in I$, $y\in J$, $x<y$.
\end{definition}

This is the standard lexicographic order if $T$ is increasing on
each element of $P$.

\begin{exercise}
Show that $\iota:(I,\leq)\to(\Sigma_+,\preceq_T)$ is order
preserving.
\end{exercise}

\begin{proposition}\label{prop-iti-admi}
Let $(T,P)$ be a piecewise monotone map. $A\in P^\NN$ belongs to
$\Sigma_+(T,P)$ iff
 \begin{equation}\label{eq-iti-admi}
   \forall n\geq0\quad \iota(e_n) \preceq_T \sigma^n A \preceq_T
   \iota(f_n)
 \end{equation}
where $e_n=c_i+$ and $f_n=c_{i+1}-$ if $(c_i,c_{i+1})=A_n$.
\end{proposition}

\begin{exercise}
Prove the above.
\end{exercise}

It is remarkable that the symbolic dynamics of piecewise monotone
maps admits such a simple description, even though it is not a
finite one as for subshifts of finite type (but this is impossible
in general as, for instance, their topological entropy can take
uncountably many values --- see \ref{exo-beta-h}). We shall see
that these subshifts are however very close to being of finite
types and share many of their properties. The notion of
\emph{quasi-finite type} described in Section \ref{sec-conclusion}
put this observation to work.

\begin{remark}
(\ref{eq-iti-admi}) is satisfied in particular by the kneading
invariants themselves. This can be formulated abstractly as a
property of a finite alphabet $\mathcal A$ where each element is
given a sign and two sequences in $\mathcal A^\NN$. The abstract
version of (\ref{eq-iti-admi}) is then necessary and sufficient to
ensure that these data are realized (up to an obvious
identification) as the kneading invariants of some piecewise
monotone map. See \cite{MS}.
\end{remark}

\subsection{Misiurewicz and Szlenk's lap numbers}

This formula says that the topological entropy of a piecewise
monotone map is combinatorial:

\begin{theorem}[Misiurewicz-Szlenk]
Let $f:[0,1]\to[0,1]$ be a piecewise monotone and piecewise
continous map. Let $P$ be its natural partition. Consider the
partition $P^n$ into $P,n$-cylinders. The topological entropy of
$f$ can be computed as:
 \begin{equation}\label{eq-Misiu-Slenk}
   h_\top(f) = \lim_{n\to\infty} \frac1n \log \#P^n
 \end{equation}
$\# P^n$ is called the \new{$n$-lap number} of $f$.
\end{theorem}

Using the Hofbauer diagram below one should be able to solve the
following:

\begin{problem}
The rate of convergence in the above formula is exponential.
\end{problem}

One can give a very easy ergodic proof of the Misiurewicz-Slenk
formula:

Recognize in eq. (\ref{eq-Misiu-Slenk}) the topological entropy of
the symbolic dynamics $\Sigma_+(f,P)$. The equality
$h_\top(\sigma|\Sigma_+(f,P))=h_\top(f)$ will follow rather easily
from the variational principle. Indeed, according to this
principle it is enough to see that one can identify invariant and
ergodic probability measures of non-zero entropy of both systems
in an entropy preserving way. Conclude by applying Exercise
\ref{exo-sd-pmm}.

\begin{exercise}\label{exo-beta-h}
Prove that if $\beta>1$ and $T:x\mapsto \beta x \mod 1$ on
$[0,1]$, then $h_\top(T)=\log\beta$ (you may consider first the
case $\beta\in\NN$). \emph{Hint:} Use Ruelle's inequality.
\end{exercise}

\subsection{Misiurewicz's Horseshoes}

\begin{definition}
A \new{horseshoe} is a collection of pairwise disjoint compact
intervals $J_1,\dots, J_N$ and an integer $T\geq1$ such that
$f^T(J_i)$ contains a neighborhood of $J_1\cup\dots\cup J_N$. The
entropy of the horseshoe is $\log N/T$.

A \new{piecewise monotone horseshoe} is a horseshoe such that, in
the above notation, $f^T|J_i$ is continuous and strictly monotone.
\end{definition}

\begin{theorem}[Misiurewicz]\label{theo-horseshoes}
Let $f$ be (i) a continuous map or (ii) a piecewise monotone map
of a compact interval. If the topological entropy of $f$ is not
zero, then it is the supremum of the entropies of its horseshoes
(its piecewise monotone horseshoes if $f$ is piecewise monotone).
\end{theorem}

Complete proofs can be found in \cite{AM} and also
\cite{RuetteBook}. These are combinatorial proofs. We give one
using ergodic theory and techniques involved in the \emph{spectral
decompositon} of A.M. Blokh \cite{BlokhDecomposition}. We consider
only the case $f\in C^0(I)$.

\begin{exercise}
Prove Misiurewicz's theorem for piecewise monotone maps either by
adapting the proof below or by using Hofbauer's Markov diagram
(see below).
\end{exercise}

\subsubsection*{Ergodic construction of Misiurewicz horseshoes}

We have to build a horseshoe of entropy at least $H$ for every
$H<h_\top(f)$. Fix such a $H$. By the variational principle, we
can find an ergodic and invariant probability measure $\mu$ such
that $h(f,\mu)>H$. By Katok's entropy formula, we can find
$\eps_0>0$ so that $h(f,\mu,\eps_0)>H$.

\step1{A Cycle}

Let $\supp'\mu$ be the support of $\mu$ minus the countably many
points $x$ such that $\mu([x,x+t])\mu([x-t,x])=0$ for $t>0$ small
enough. In the rest of the construction, $x$ will be a point of
$\supp'\mu$.

If $I$ is an open interval, define $F(I)$ as the interior of
$f(I)$. For $r>0$ and $x\in[0,1]$, define:
 $$
    V_r(x) := \bigcup_{k\geq 0} F^k(B(x,r)) \text{ and }
    C_r(x) \text{ is the c.c. of $V_r(x)$ containing $x$}.
 $$
As $x\in\supp'\mu$, there exists $k\geq1$ such that
$F^k(B(x,r))\cap B(x,r)$ has positive measure. Hence $V_r(x)$ has
finitely many, say $n(x,r)$, connected components.

\begin{exercise}
Check that $F$ permutes these connected components and that $x\in
V_r(x)$.
\end{exercise}

We claim that
 \begin{equation}\label{eq-cycle}
   \sup_{r>0} n(x,r)<\infty.
 \end{equation}
An entropy argument shows that the proportion $\xi(x,r)$ of $0\leq
k<n(x,r)$ such that $\diam(F^kC_r(x))>\eps_0$ satisfies:
$\xi(x,r)\geq H/\log(\eps_0^{-1}+1)$. Indeed, $r(\eps_0,n,C_r(x))
\leq (\eps_0^{-1}+1)^{\xi(x,r) n}$ so $H<h(f,\mu,\eps_0)\leq
\xi(x,r) \log(\eps_0^{-1}+1)$. But, considering the lengths of the
pairwise disjoint $F^k(C_r(x))$, we have:
$(1-\xi(x,r))n(x,r)\eps_0\leq 1$, proving eq. (\ref{eq-cycle}).
This implies that $\bigcap_{r>0} C_r(x)=:C(x)$ is a non-trivial
interval defining a cycle $V(x)$ of period $p(x)$ with
$F^{p(x)}(C(x))=C(x)$ (if $F^{p(x)}(C(x))\subsetneq C(x)$, then
taking $r>0$ so small that $F^{p(x)}(C_r(x))\cup B(x,r)\subsetneq
C(x)$ we would get $C_r(x)\subsetneq C(x)$, a contradiction).

If $y$ is another point of $\supp'\mu$, we claim that $V(y)=V(x)$.
Indeed, observe that $\mu(V(x))=1$ by ergodicity, so
$\supp\mu\subset \overline{V(x)}$. As $y\in\supp'\mu$, $y\in
V(x)$. Hence $B(y,r)\subset V(x)$ for $r>0$ small enough, so
$V(y)\subset V(x)$. By symmetry, $V(y)=V(x)$. Henceforth, we drop
the dependence on $x$, simplifying $V(x),C(x),p(x)$ to $C,V,p$.

\step2{Strong Mixing}

Call a subinterval $I\subset [0,1]$ such that $\bigcup_{k\geq0}
F^k(I)\subsetneq V$ \new{negligible}. Introduce the following
relation on $[0,1]$: $x\sim y$ iff $[x,y]$ is negligible. Observe
that $x\sim y$ implies $f(x)\sim f(y)$ as $[f(x),f(y)]\subset
f([x,y])$ and that the equivalence classes of $\sim$ are compact
subintervals.

\begin{exercise}
Show that $f:[0,1]\to[0,1]$ induces a continuous self-map of
$[0,1]/\sim$ and that, up to a homeomorphism, this new map is a
continuous map $f'$ of $[0,1]$ which admits an invariant and
ergodic probability measure $\mu'$ satisfying
$h(f',\mu',\eps_0')>H$ for some $\eps_0'>0$.
\end{exercise}

By this exercise one can assume that no non-trivial subinterval of
$V$ is negligible. Hence for every non-trivial subinterval
$J\subset C$, $\bigcup_{k\geq0} F^kJ=V$.

Observe that $f^p(C)=C$ implies the existence of a fixed point
$f^p(x_0)=x_0$ in the closure $\bar C$. In particular, there is an
integer $n_0>0$ such that $f^{n_0}J\ni x_0$ by the following:

\begin{exercise}
Show that either one can find $x_0\in C$ (in which case we set
$y_0:=x_0$), or there exists $y_0\in V$ such that $f(y_0)=x_0$.
\end{exercise}

Observe that $f^{n_0}J$ cannot be negligible, as it would imply
that all connected components of $\bigcup_{k=0}^{n_0-1}f^kJ$ must
also be negligible, but this union must have full $\mu$-measure, a
contradiction.

If $x_0\in\partial C$, let $z_0\in C$ such that
$\mu([x_0,z_0])>0$. Otherwise, let $w_0\in\partial C$ such that
$\mu([x_0,w_0])>0$ and $z_0\in(x_0,w_0)$ such that
$\mu([x_0,z_0])\mu([z_0,w_0])>0$. In both cases, set
$K_0:=[x_0,z_0]$.

As $f^{n_0}J$ is not negligible, there exists $n_1>n_0$ such that
$f^{n_1}J\ni z_0$, so $f^{n_1}J\supset K_0$.

\step3{Conclusion}

Divide $V$ into subintervals $I_1,\dots,I_L$, each of length
$\eps_0/2$. By Step 2, for each $i=1,\dots,L$, a positive integer
$t_i$ such that $f^{t_i}(I_i)\supset K_0$. Let $t_*:=\max t_i$.

By Katok's entropy formula, $\mu(K_0)>0$ implies that, for all
large $n$, there exist $x_1<x_2<\dots<x_N$ points in $K_0$ which
are $(\eps_0,n)$-separated with $N\geq 2ne^{t_*H}e^{nH}$: for
every $i=1,\dots,N-1$, there exists $n_i<n$ such that
$|f^{n_i}x_i-f^{n_i}x_{i+1}|\geq\eps_0$. Hence,
$f^{n_i}([x_i,x_{i+1}])\supset I_j$ for some $j=j_i\in\{1,2,\dots,
L\}$. Thus, $f^{n_i+t_j}([x_i,x_{i+1}])\supset K_0$.

The integers $n_i+t_{j_i}$ belong to $[0,n+t_*]$. Hence, one can
select $i_1,\dots,i_M$ from $\{1,\dots,N\}$ with $M=N/(n+t_*)\geq
e^{(n+t_*)H}$, $j(i_1)=\dots=j(i_M)$ and $T:=n_{i_1}+t_{j(i_1)}$
such that:
 \begin{itemize}
  \item $J_\ell:=I_{i_\ell}\subset\subset K_0$ for $\ell=1,\dots,M$;
  \item $J_\ell\cap J_{\ell'}=\emptyset$ if $\ell\ne\ell'$;
  \item $f^T(J_i)\supset K_0$
 \end{itemize}
That is, $f$ has a horseshoe with entropy at least $H$, proving
Misiurewicz theorem.

\subsection{Entropy and length}

Entropy is reflected in the growth rate of volume. The situation
is especially simple on the interval \cite{MP,Yomdin}.

\begin{proposition}[Misiurewicz-Przyticky, Yomdin]
For $f\in C^1(I)$ define the growth of length under iterations of
$f$ to be:
 $$
   \ell(f) := \limsup_{n\to\infty} \frac1n\log \int_I
   |(f^n)'(x)|\, dx.
 $$
Then $\ell(f)\geq h_\top(f)$. For $f\in\PMM(I)\cap C^1(I)$ or
$f\in C^\infty(I)$, this is in fact an equality:
 \begin{equation}\label{eq-htop-len}
    h_\top(f) = \ell(f).
 \end{equation}
\end{proposition}

The above remain true for $f\in\PMM(I)$, provided that $\ell(f)$
is suitably redefined using the variation of $f^n$ (minus that
occuring at its discontinuities). The proof is the same for
$f\in\PMM(I)\cap C^0(I)$ and slightly more delicate for the
general case.

This formula has generalizations to $C^\infty$-smooth self-maps in
arbitrary dimension -- see \cite{Newhouse1} and also
\cite{Kozlovski}. One can compare it with the following
consequence of Rokhlin formula (Theorem \ref{theo-rokhlin}). Let
$f$ be a piecewise monotone, piecewise $C^1$ map of the interval.
If $\mu$ is an invariant probability measure with absolutely
continuous ergodic components, then:
 $$
    h(f,\mu) = \int \limsup_{n\to\infty}(1/n)\log|(f^n)'(x)| \, d\mu
 $$
This has been generalized to arbitrary invariant probability
measures absolutely continuous (this is Pesin's formula) or not
and self-maps on manifolds of arbitrary dimension by taking into
account the relevant \emph{dimensions} (see
\cite{LSY-D,Hofbauer-D}) and considering the \emph{Lyapunov
exponents} instead of the single derivatives (see
\cite{KatokStrelcyn,LY}).

\begin{exercise}
Check (\ref{eq-htop-len}) in the case where $f$ maps each of its
piece to the full interval (up to its endpoints).
\end{exercise}

\begin{demo}
We first prove $\ell(f)\geq h_\top(f)$ for all $C^1$ maps. Let
$\eps>0$. For all large $n$, there exists $x_1<\dots <x_N$ a
$(\eps,n)$-separated subset with $N\geq
ne^{n(h_\top(f)-\eps)}/\eps$. As in the proof of Theorem
\ref{theo-horseshoes}, we can find a subset $y_1<\dots<y_M$ with
$M\geq N/n$ and $0\leq m\leq n$ such that
$|f^ky_i-f^ky_{i+1}|\geq\eps$. Therefore,
 $$
   \frac1n\log \int_I
   |(f^m)'(x)|\, dx \geq e^{m(h_\top(f)-\eps)}
 $$
which implies the claimed inequality.

We prove the converse inequality. The $C^\infty$ case follows from
Yomdin's theory - Theorem \ref{theo-yomdin} - we leave its details
to the diligent reader. For $f\in\PMM(I)\cap C^1(I)$, we observe:
 $$
  \int_I
   |(f^n)'(x)|\, dx = \sum_{A\in P^n} |f^n(A)| \leq \#P^n
 $$
so that the Misiurewicz-Szlenk formula is enough to conclude.
\end{demo}

\begin{remark}
Misiurewicz and Przytycki \cite{MP} have shown that, in any
dimensions, the logarithm of the degree is a lower bound for the
entropy for $C^1$ maps (but not $C^0$ maps). This is a (rather)
special case of the Shub entropy conjecture mentioned before.
\end{remark}

\section{Entropy as a function}\label{sec-function-entropy}

\subsection{Local entropy}
Recall that local entropy quantifies the complexity at arbitrarily
small scales. The situation is not so bad for "reasonable"
interval maps:

\begin{proposition}\label{prop-hloc0}
For $f\in\PMM(I)$ or $f\in C^\infty(I)$, $h_\loc(f)=0$.
\end{proposition}

\begin{exercise}
If $K\subset I$ is such that $f^k|K$ is monotone for $0\leq k<n$,
then $r(\eps,n,K)\leq n(|I|/\eps+1)+1$.
\end{exercise}

\begin{demo}
We have already explained how the assertion for $f\in C^\infty(I)$
follows from Yomdin's theory. Assume $f\in\PMM(I)$.

Consider first $f\in\PMM(I)$. The exercise below implies that
$r(\eps,n,J) \leq Cn \#\{A\in P^n:A\cap J\ne\emptyset\}$ where $P$
is the natural partition of $f$. Now, if $K\subset B(\eps,n,x)$,
then each $f^kK$, $0\leq k<n$, meets at most two elements of $P$,
hence $\#\{A\in P^n:A\cap J\ne\emptyset\}\leq 2^n$ so that
$h_\loc(f)\leq 2$. Applying the same to $f^N$ we get, by Exercise
\ref{exo-hloc-iter} $h_\loc(f)=(1/N) h_\loc(f^N)\leq \log 2/N$.
Letting $N\to\infty$, we get the claim for piecewise monotone
maps.
\end{demo}

\subsection{Measured entropy}

Recall that a map $\phi:X\to\mathbb R$ is \emph{lower
semi-continuous} at $x\in X$, if, for $\eps>0$, there exists a
neighborhood of $x$ on which $\phi$ is at least $\phi(x)-\eps$:
that is, $\phi$ cannot "collapse". \emph{Upper semi-continuity} of
$\phi$ is lower semi-continuity for $-\phi$: it means that $\phi$
cannot "explode". In particular an upper semi-continuous function
over a non-empty compact subset achieves its maximum.

The preceding result on the local entropy yields the:

\begin{corollary}
For any $f\in\PMM(I)$ or $f\in C^\infty(I)$, the measured entropy
is upper semi-continuous on the compact of invariant probability
measure and therefore $f$ admits measures of maximum entropy.
\end{corollary}

Note that such an existence proof says nothing about the
multiplicity or the structure of these measures. These two points
will require a rather complete analysis of the one-dimensional
dynamics through Hofbauer's Markov diagram below.

\begin{fact}
If $f\in C^0(I)$ or $f\in \PMM(I)$ and $h(f,\mu_0)>0$, then
$\mu\mapsto h(f,\mu)$ is not lower semi-continuous at $\mu=\mu_0$.
\end{fact}

This follows from the density of periodic measures (invariant and
ergodic probability measures defined by periodic orbits) among all
invariant measures (which holds because of specification
properties enjoyed by these maps (see
\cite{BlokhDecomposition,BuzziSpec,HofbauerSpec})).

\begin{fact}
If $f\in C^\infty(I)$ or $f\in\PMM(I)$, then $\mu\mapsto h(f,\mu)$
is upper semi-continuous.

For any $r<\infty$ there are examples of $f\in C^r(I)$ with
invariant probability measures $\mu_0$ at which the entropy fails
to be upper semicontinuous.
\end{fact}

The first fact follows from $h_\loc(f)=0$ (see Proposition
\ref{prop-hloc0}). For $r<\infty$, examples on the interval are
given in \cite{BuzziSIM, RuetteMax}. {\it Note:} In these
examples, $h(f,\mu_0)=0$.

\begin{remark}
S. Ruette \cite{RuetteBorel} observed that this gives examples of
topologically mixing, smooth interval maps with the same nonzero
topological entropy which are not  Borel conjugate.
\end{remark}

\begin{remark}
M. Rychlik \cite{Rychlik1} proved that the measured entropy of the
unique absolutely continuous invariant probability measure of a
family of piecewise $C^2$, piecewise expanding map has modulus of
continuity $x\log1/x$.
\end{remark}

The following should be tractable by using the perturbative
results of Keller and Liverani \cite{KellerLiverani} about the
transfer operator.

\begin{problem}
Study the regularity of the measured entropy for families of
piecewise expanding, piecewise monotone maps.
\end{problem}

\subsection{Topological entropy}

We consider the dependance of the topological entropy on the map.

\subsubsection{Continuity}

\begin{theorem}
The map $f\mapsto h_\top(f)$ is lower semi-continuous both over
$C^0([0,1])$ and over $\PMM([0,1])$.

This map is upper semi-continuous (and therefore continuous) over
$C^\infty([0,1])$ and over $\PMM([0,1])$ but not over $C^r([0,1])$
for any finite $r$.
\end{theorem}

\begin{remark}
In our choice of topology over $\PMM(I)$, two maps can be close
only if they have the same number of pieces. This is not always
the chosen definition but it is necessary to get the upper
semi-continuity.
\end{remark}

It is easy to find a counter-example to the upper semi-continuity
in the $C^0(I)$ topology even among piecewise monotone maps:

\begin{exercise}
Let $f_t(x)=(1-t)(1-2t|x|/(1-t))$ if $|x|<1-t$ and $f_t(x)=0$
otherwise. Show that $h_\top(f_t)=\log(2t)$ pour $t<1$ and
$h_\top(f_1)=0$.

Can you find such an example among polynomials?
\end{exercise}

The basic phenomenon here is that a periodic point becomes a cycle
of intervals containing a horseshoe. This in fact the only
obstruction -- see \cite[Prop. 4.5.3]{AM}. We refer to \cite[Th.
9.1]{MS} or \cite[Cor. 4.5.5]{AM} for two rather different proofs:
the original one, due to Milnor and Thurston, which involves the
"kneading determinant" which is a sort of zeta function, and a
combinatorial one, closer to the spirit of this paper. We indicate
below the outline of the proof in the smooth case.

 \medbreak

The proof of the above theorem will use the:

\begin{exercise}\label{exo-Cfn-C0}
Show that for every $\eps>0$, $C(g^n)\subset B(C(f^n),\eps)$ for
all $g$ close enough to $f\in\PMM(I)$.
\end{exercise}

\begin{exercise}
Assume that for all $\alpha>0$, there exist $\eps_0>0$ and a
continuous function $C_\alpha:(0,\infty)\to(0,\infty)$ such that,
for all $g$ in a neighborhood of $f$, all $\delta>0$:
 \begin{equation}\label{eq-unif-hloc}
   \forall n\geq 1\; r_g(\delta,n,B_g(\eps_0,n,x)) \leq C_\alpha(\delta)
   e^{\alpha n}.
 \end{equation}
Show that $g\mapsto h_\top(g)$ is then upper semi-continuous at
$g=f$.
\end{exercise}

\begin{demo}
The lower semi-continuity follows from Theorem
\ref{theo-horseshoes}: recall that a horseshoe is defined by
finitely many conditions like (i) $f^n|[a,b]$ continuous; (ii)
$\INT f^n([a,b])\supset[c,d]$ and each of these is equivalent to
$f^n|[a,b]$ continuous and $f^n(a')<c$ and $f^n(b')>d$ for some
$a',b'\in [a,b]$. This is obviously an open condition in $C^0(I)$.
It is also open in $\PMM(I)$ by Exercise \ref{exo-Cfn-C0}.

The upper semi-continuity of the topological entropy in the
$C^\infty$ smooth case follows from Yomdin's theory which gives,
not only that $h_\loc(f)=0$, but the very strong uniformity stated
as the assumption in (\ref{eq-unif-hloc}). By the above exercise,
this concludes the proof of the proposition.

We refer the reader to the references quoted above for the proof
in the piecewise monotone case.
\end{demo}

Lower semi-continuity also holds for surface diffeomorphisms of
class $C^{1+\eps}$ according to a classical result \cite{Katok0}
of Katok which shows the existence of horesehoes. This is false
for surface homeomorphisms as can be deduced from examples of
Rees. The natural case of $C^1$ diffeomorphisms is still open. It
fails however in higher dimensions, even in the analytic category.

\begin{exercise}
Find a family of maps $F_\lambda:[0,1]^2\to[0,1]^2$,
$\lambda\in[0,1]$ with $F_\lambda(t,x)=(\lambda t,f_t(x))$ such
that $\lambda=1$ is a discontinuity point for $\lambda\mapsto
h_\top(F_\lambda)$.
\end{exercise}

The following question seems not to have been studied though the
characterization of entropy as an eigenvalue of a transfer
operator should allow the application of techniques like
\cite{KellerLiverani}.

\begin{question}
Study the modulus of continuity of topological entropy in families
of piecewise monotone maps.
\end{question}

\subsubsection{Monotonicity}

The following fact was proved by Sullivan, Milnor, Douady and
Hubbard:

\begin{theorem}
Let $Q_t(x)=4tx(x-x)$ be the quadratic family $Q_t:[0,1]\to[0,1]$.
Then $t\mapsto h_\top(Q_t)$ is non-decreasing for $t\in[0,1]$.
\end{theorem}

\begin{remark}
This sounds very natural. One must be careful however that the
similar statements involving other families of the type $t\mapsto
tQ(x)$ for a smooth map $Q:[0,1]\to[0,1]$ satisfying $Q'(x)=0$ iff
$x=1/2$ for $x\in[0,1]$. See \cite{BruinNonMono}. See
\cite{Misiurewicz-Visinescu,Galeeva,MiTr} for further results.
\end{remark}

Most of the proofs rely on complex dynamics (e.g., Teichmuller
theory or quasiconformal maps). Let us describe the strategy of
the remarkable "real"\footnote{With the following restriction: as
remarked by Tsujii in his paper, the map $R$ below is a local
version of Thurston map on the Teichmuller space - see
\cite{Levin}.} proof found by Tsujii \cite{TsujiiMonotonocity}.

Recall how the kneading invariants of a piecewise monotonic map
define its symbolic dynamics (Proposition \ref{prop-iti-admi}) and
therefore its topological entropy. Observe that $Q_t(0)=Q_t(1)=0$
and $Q_t$ is continuous, hence there is only one non-trivial
kneading invariant, $\kappa(Q_t):=\iota(Q_t(1/2))$ and that, with
respect to the lexicographic order defined by any of the $Q_t$,
the symbolic dynamics as a subset and therefore the topological
entropy as a number are non-decreasing functions of $\kappa(Q_t)$.

Proposition \ref{prop-iti-admi} implies that the map that
associates to such a kneading invariant the topological entropy of
the corresponding unimodal map is nondecreasing. Hence
monotonicity of $t\mapsto \kappa(Q_t)$ (for the above order on
$\{-1,0,1\}^\NN$ will imply the monotonicity of the entropy. Hence
the above theorem is reduced to showing the following property:
 \begin{equation}\label{eq-mono}
   \text{if }n:=\min\{k\geq1:Q_t^k(0)=0\} \text{ then, }
      \frac{\partial}{\partial t} Q_t^n(0) \text{ and }
      \frac{\partial}{\partial x} Q_t^n(Q_t(0))
 \end{equation}

Tsujii proves this last property by considering the following
\new{Ruelle transfer operator} (see \cite{BaladiBook} for this very important notion):
 $$
     R(\Psi)(x) = \sum_{y\in Q^{-1}(x)} \frac{\Psi(y)}{(\partial
     Q_t(y)/\partial y)^2}
 $$
acting on the space $E$ linearly generated by $1/(z-Q^k(0))$,
$k\geq 0$ (there are exactly $n$ such functions because of the
assumption on $t$). A computation shows that
 $$
 \det(\operatorname{Id}_E -R|E)=\frac{\partial}{\partial t}
Q_t^n(0)/\frac{\partial}{\partial x} Q_t^n(Q_t(0))
 $$
Now, Tsujii proves that $R|E$ is a contraction. This proves that
$\phi(z):= \det(\operatorname{Id}_E -zR|E)\ne 0$ for all
$|z|<\rho$ for some $\rho>1$. Clearly $\phi(0)=1$ and $\phi(z)$ is
real for all real $z$. Hence, $\phi(1)>0$, concluding Tsujii's
proof.

\section{The Markov Diagram}\label{sec-diagram}

The main tool for the detailed analysis of the complexity of
interval maps with nonzero entropy is Hofbauer's \new{Markov
diagram}, which reduces much of their study to that of a countable
state Markov shift with very good properties. It is inspired from
a classical construction in the theory of subshifts of finite type
and sofic subshifts: \emph{the canonical extension}, which builds
for any sofic subshift a topological extension which is a subshift
of finite type with the same entropy --- see \cite{LM}.

\subsection{Abstract construction}

Let us describe this construction, or rather a variant introduced
in \cite{BuzziSIM} which allows slightly stronger statements (and,
more importantly, is necessary for generalizations, see remark
\ref{rem-original-MD} below).

The Hofbauer construction applies to any subshift. Let
$\Sigma_+\subset\mathcal A^\NN$ be a one-sided subshift (i.e., a
closed, shift-invariant subset of $\mathcal A^\NN$, not
necessarily of finite type). Let $\Sigma\subset\mathcal A^\ZZ$ be
its natural extension:
 $$
   \Sigma:=\{A\in\mathcal A^\ZZ:\forall p\in\ZZ\;
   A_pA_{p+1}\dots\in \Sigma_+\}
 $$
with the action of the (invertible) shift $\sigma: A\mapsto
(A_{n+1})_{n\in\ZZ}$

\begin{definition}
A \new{word} is a map from a finite integer interval
$\{a,a+1,\dots,b\}$ to $\mathcal A$, up to an integer translation
of the indices. The
\new{follower set} of a word $A_{-n}\dots A_0$ is:
 $$
   \fol(A_{-n}\dots A_0) = \{B_0B_1\dots \in\Sigma_+:\exists B\in\Sigma\;
   B_{-n}\dots B_0=A_{-n}\dots A_0\}.
 $$
A \new{minimal word} of $\Sigma$ (or $\Sigma_+$) is $A_{-n}\dots
A_0$ such that:
 $$
   \fol(A_{-n}\dots A_0) \subsetneq \fol(A_{-n+1}\dots A_0).
 $$
The \new{minimal form} of $A_{-n}\dots A_0$ is
 $$
   \min(A_{-n}\dots A_0):=A_{-k}\dots A_0
 $$
where $k\leq n$ is maximum such that $A_{-k}\dots A_0$ is minimal.
\end{definition}

\begin{definition}
The \new{Markov diagram} $\MD$ of a subshift $\Sigma$ is the
oriented graph whose vertices are the minimal words of $\Sigma$
and whose arrows are:
 $$
   A_{-n}\dots A_0 \to B_{-m}\dots B_0 \iff
     B_{-m}\dots B_0 = \min(A_{-n}\dots A_0B_0).
 $$
The corresponding Markov shift is:
 $$
   \HS:=\{\alpha\in\MD^\ZZ:\forall p\in\ZZ\;
     \alpha_p\to\alpha_{p+1} \text{ on }\MD\},
 $$
with the left shift $\sigma$ may be called the \new{Hofbauer
shift}. There is a natural, continuous projection:
 $$
   \hat\pi: \alpha\in\HS\longmapsto A\in\Sigma
 $$
with $A_n$ the last symbol of the word $\alpha_n$ for each
$n\in\ZZ$.
\end{definition}

\begin{remark}\label{rem-original-MD}
The vertices of Hofbauer's original diagram were simply the
follower sets, as in the classical construction for sofic
subshifts \cite{LM}. It can also be used beyond piecewise
monotonic maps \cite{BuzziMDD} but it leads to a less tight
relationship of the dynamics. In particular, periodic points may
have several lifts\footnote{This creates only a finite number of
extra periodic orbits on the interval but can be a serious issue
for generalizations}.
\end{remark}

\begin{exercise}
Show that a subshift has a finite Markov diagram iff it is a
subshift of finite type. \emph{Note:} Hofbauer's Markov diagram is
finite if and only if the subshift is sofic.
\end{exercise}

\begin{exercise}
Compute the Markov diagram of the even shift, that is,
$\Sigma\subset\{0,1\}^\ZZ$ obtained by excluding the words
$01^{2n+1}0$, $n\geq0$.
\end{exercise}

\subsection{Eventual Markov property and partial isomorphism}

We are going to see that the Hofbauer shift $\HS$ is partially
isomorphic with $\Sigma$, the (invertible) original subshift.

\begin{definition}
Let $\Sigma$ be a subshift on a finite alphabet. $A\in\Sigma$ is
\new{eventually Markov} (or just Markov) at time $p\in\ZZ$ if there exists $N=N(x,p)$ such that:
 $$
    \forall n\geq N\; \fol(A_{p-n}\dots A_p) = \fol(A_{p-N}\dots A_p)
 $$
The \new{eventually Markov part} $\Sigma_M\subset\Sigma$ is the
set of $A\in\Sigma$ which are eventually Markov at all times
$p\in\ZZ$.
\end{definition}

\begin{remark}
$\Sigma_M$ is a topological version of the so-called "variable
length Markov chains".
\end{remark}

The following was shown by Hofbauer for piecewise monotone maps
and the original Markov diagram (building on work of Y. Takahashi)
and then generalized to arbitrary subshifts
\cite{BuzziSIM,BuzziMDD}:

\begin{theorem}[Hofbauer, Buzzi]\label{theo-partial-iso}
The natural projection from the Hofbauer shift $\HS$ to the
subshift $\Sigma$ defined by:
 $$
    \hat\pi:(\alpha_p)_{p\in\ZZ} \longmapsto (A_p)_{p\in\ZZ}
    \quad\text{with $A_p\in P$ the last symbol of $\alpha_p$}
 $$
is well-defined and is a Borel isomorphism from $\HS$ to
$\Sigma\setminus\Sigma_+$.
\end{theorem}

\begin{exercise}
Show that if $\Sigma$ is the symbolic dynamics defined by an
irrational rotation $T$ with $P$ a partition into two disjoint
intervals, then $\Sigma_M$ is empty.
\end{exercise}

\begin{lemma}
Let $\alpha_0\to\dots\alpha_n$ be a finite path on $\MD$. Let
$B_{-k}\dots B_0=\alpha_0$ and $A=\hat\pi(\alpha)\in\Sigma$. Then,
 $$
   \alpha_n = \min(B_{-k}\dots B_0A_1A_2\dots A_n).
 $$
\end{lemma}

\begin{exercise}
Prove the lemma.
\end{exercise}

\begin{demof}{Theorem \ref{theo-partial-iso}}
The following measurable map is a candidate to be a partial
inverse to $\hat\pi$:
 $$
   \iota: A\in \Sigma_M\longmapsto \alpha\in\HS
   \quad\text{with }\alpha_p=\lim_{n\to\infty}\min(A_{p-n}\dots
   A_p)
 $$
The above limit is to be understood in the most trivial sense: it
is equal to $A_{p-N}\dots A_p$ if this is the value of
$\min(A_{p-n}\dots A_p)$ for all large $n$ (i.e., use the discrete
topology on the countable set of finite words). $A\in\Sigma_M$ is
equivalent to the existence of this "limit" for all $p\in\ZZ$.

\step1{$\hat\pi:\HS\to\Sigma_M$ is well-defined}

Let $\alpha\in\HS$. Let $A:=\hat\pi(\alpha)\in\mathcal A^\ZZ$. We
first show that for each $p\in\ZZ$ and $n\geq0$, $A_p\dots
A_{p+n}$ is a word of $\Sigma_+(T,P)$. But the lemma applied to
$\alpha_p\to\alpha_{p+1}\dots$ implies that it is a suffix of such
a word, which gives the property so $A\in\Sigma$.

Let us check that $A\in \Sigma_M$. For any $n\geq 0$,
$\alpha_{-n}=D^n_{-k(n)}\dots D^n_0$ with $D^n_0=A_{-n}$. Let
$m:=k(0)$. By the Lemma applied to
$\alpha_{-n}\to\dots\to\alpha_0$, we have:
 \begin{equation}\label{eq-iota-reduction}
    \min(D^n_{-k(n)}\dots D^n_{-1}A_{-n}\dots A_0) = D^0_{-m}\dots D^0_0
 \end{equation}
But this implies that
 $$
  \forall n\geq m\; \min(A_{-n}\dots A_0) = D^0_{-m}\dots D^0_0
 $$
so that $A$ is Markov at time $0$. The same applies to any time,
hence $A\in\Sigma_M$.

\step2{$\iota(\Sigma_M)\subset\HS$}

Let $A\in\Sigma_M$ and $\alpha:=\iota(A)\in\MD^\ZZ$. Let
$p\in\ZZ$. We have, for $n$ large enough,
 $$
   \alpha_p = \min(A_{p-n}\dots A_p) \text{ and }
   \alpha_{p+1} = \min(A_{p+1-n}\dots A_{p+1})
 $$
Hence $\alpha_p=A_{p-\ell}\dots
A_p=\min(A_{p-\ell-1}A_{p-\ell}\dots A_p)$, so:
 $$
   \sigma([A_{p-\ell-1}]) \supset [A_{p-\ell}\dots A_p]
 $$
which implies:
 $$
  \sigma([A_{p-\ell-1}]) \supset [A_{p-\ell}\dots A_{p+1}]
 $$
giving that $\min(A_{p-\ell-1}\dots A_{p+1})= \min(A_{p-\ell}\dots
A_{p+1})$. By induction:
 $$
   \alpha_{p+1}=\min(A_{p-n}\dots A_{p+1})=\min(A_{p-\ell}\dots
A_{p+1})
 $$
We have shown:
 $$
   \alpha_p=A_{p-\ell}\dots A_p \text{ and }
   \alpha_{p+1} = \min(A_{p-\ell}\dots
A_{p+1})
 $$
and this is the definition of $\alpha_p\to\alpha_{p+1}$. Hence
$\iota(\Sigma_M)\subset \HS$.

\step3{$\hat\pi:\HS\to\Sigma_M$ is a bijection}

By the previous steps, $\iota\circ\hat\pi(\alpha)$ is well-defined
for $\alpha\in\HS$. Eq. (\ref{eq-iota-reduction}) implies that
 $$
   \alpha_0 = D^0_{-m}\dots D^0_0 = A_{-m}\dots A_0
   = \min(A_{-n}\dots A_0) \quad \forall n\geq m.
 $$
Hence, $\iota\circ\hat\pi:\HS\to\HS$ is the identity.

The previous steps also show that $\hat\pi\circ\iota$ is
well-defined over $\Sigma_M$ and it is obviously the identity
there. This concludes the proof of the theorem.
\end{demof}

\begin{remark}
G. Keller introduced another point of view in \cite{KellerLifting}
by observing that the one-sided Markov shift defined by the Markov
diagram contains a (non necessarily invariant) copy of $\Sigma$ so
that one can lift any probability measure of $\Sigma$ there and
try and make it invariant by pushing it forward. It leads to other
type of lifting theorems. See Remark \ref{rem-exp-shad} below for
a (quick) comparison between the two approaches.
\end{remark}

\subsection{Markov diagrams for piecewise monotone maps}

Hofbauer discovered that the Markov diagram of a piecewise
monotone map is very special. Recall how the symbolic dynamics is
defined by the finite set of
\new{kneading invariants} $K(T,P)$ (see Section \ref{sec-kneading}
above).

\begin{theorem}[Hofbauer]\label{theo-PMM-MD}
Let $T:[0,1]\to[0,1]$ be a piecewise monotone map with natural
partition $P$. Let $\Sigma_+(T,P)$ be its symbolic dynamics.

All the vertices of the Markov diagram are obtained from the
kneading invariants:
 $$
   \MD =\{A_0\dots A_n: 0\leq n<N(A)\text{ and }A\in K(T,P)\}
 $$
where $N(A):=\inf\{n\geq0:\min(A_0\dots A_n)\ne A_0\dots A_n\}$
(usually infinite) and the arrows are exactly:
 \begin{enumerate}[(i)]
  \item $A_0\dots A_n\to\min(A_0\dots A_{n+1})$;
  \item $A_0\dots A_{S(A,i+1)-1}\to \min(A_{S(A,i)}\dots
  A_{S(A,i+1)})$ if $S(A,i+1)<\infty$;
  \item $A_0\dots A_{S(A,i+1)}\to Q$ for $Q\in P(A,i)$  if $S(A,i+1)<\infty$.
 \end{enumerate}
where $P(A,i)\subset P$ and
$S:K(T,P)\times\NN\to\{1,2,\dots,\infty\}$ satisfies for each
$i\geq0$, $S(A,i+1)=S(A,i)+S(B,j)$ or $\infty$ (in the finite case
$(B,j)\in K(T,P)\times\NN$ is determined by $B_0\dots
B_{S(B,j)}=A_{S(A,i)}\dots
  A_{S(A,i+1)}$).
\end{theorem}

\begin{exercise}
For $\beta>1$, compute the Markov diagram of the
\new{$\beta$-transformation}: $T_\beta:[0,1]\to[0,1]$ defined by
$T_\beta(x)=\beta x -[\beta x]\in[0,1)$ ($[\cdot]$ is the integer
part). \emph{Hint:} Introduce
$\{R_1<R_2<\dots\}:=\{r\geq1:\iota(1^-)_n\ne (0,\beta^{-1})\}$ and
show that, up to finitely many vertices, it can be represented as
a "linear graph" $0\to1\to 2\to\dots$ with additional "backwards"
arrows that you will specify.

Show that the Markov diagram of $T_\beta$ can be finite only if
$\beta$ is algebraic.
\end{exercise}

\begin{proposition}[Hofbauer]
Let $T:[0,1]\to[0,1]$ be \emph{unimodal}, that is, $T$ is
continuous and there exists $c\in(0,1)$ such that $T|[0,c]$ and
$T|[c,1]$ are continuous and strictly monotone, normalized by
$T(0)=T(1)=0$. Assume that $c<T(c)<1$ (otherwise the situation is
easy to analyze). Let $A=\iota(c^-)$ be the unique nontrivial
kneading invariant.

The Markov diagram of $T$ has the following structure.

Its vertices are: $I_n:=\min(A_0\dots A_n)$, for $n\geq 0$,
together with a transient vertex $(c,1)$.

Its arrows are:
 \begin{itemize}
  \item $(c,1)\to I_0,I_1$;
  \item $I_n\to I_{n+1}$ for all $n\geq0$;
  \item $I_{S(n)-1}\to I_{S(Q(n+1))}$ if $Q(n+1)<\infty$
 \end{itemize}
where $S:\NN\to\NN$ and $Q:\NN\to\NN$ satisfy
$S(n+1)=S(n)+S(Q(n+1))$ if $Q(n)<\infty$, $Q(n)<n$ unless it is
infinite.
\end{proposition}

\begin{remark}
The simplicity already observed on the level of symbolic dynamics
(which gives rise to the notion of \emph{quasi finite type}
described in section \ref{sec-conclusion}) is again visible in the
Markov Diagrams of piecewise monotone maps. They are not finite
(except in very special cases) but they are "almost so". This is
made precise by the notion of \emph{strongly positive recurrence}
(see section \ref{sec-SPR}).
\end{remark}

\begin{remark}
F. Hofbauer translated the relatively simple admissibility
condition on the kneading invariants to one on the above
$Q$-function. This results allowed the construction of quadratic
map with finely tuned properties (see, e.g.,
\cite{HKOscillations}).
\end{remark}

\section{Non-Markov sequences, Shadowing and Weak Rank
One}\label{sec-shadow}

We have showed that any invertible subshift $\Sigma$ contains a
Markov part $\Sigma_M$ which is the continuous, Borel isomorphic
image of the Markov shift defined by the Markov diagram of
$\Sigma$. To use this, one needs to control the complement
$\Sigma\setminus\Sigma_M$. For the symbolic dynamics of piecewise
monotone maps (and many other kind of subshifts, see Section
\ref{sec-conclusion}) we use the following notion.

\subsection{Non-Markov Sequences}

\begin{definition}
Let $\Sigma_+\subset\mathcal A^\NN$ be a subshift over a finite
alphabet and $\Sigma$ its natural extension. Let
$S\subset\Sigma_+$. $A\in\Sigma$ is
\new{shadowed} by $S$ if there exist arbitrarily large integers
$n$ such that $A_{-n}\dots A_0$ is the beginning of a sequence in
$S$. A measure $\mu$ on $\Sigma$ is shadowed if $\mu$-a.e.
sequence is.
\end{definition}

\begin{theorem}[after Hofbauer]\label{theo-smallentropy1}
If $(T,P)$ is a piecewise monotone map, then the non-Markov part
of its invertible symbolic dynamics $\Sigma(T,P)$ carries only
zero entropy measure and therefore the partial isomorphism of
Theorem \ref{theo-partial-iso} is an entropy-conjugacy.

More precisely, except for finitely many periodic sequences, a
sequence of $\Sigma(T,P)$ is shadowed at all times by $K(T,P)$ iff
it is non-Markov at all times.

Moreover, $\hat\pi$ is a period-preserving bijection between the
periodic orbits of $\HS$ and their image in $\Sigma_+(T,P)$ except
for finitely many periodic orbits.
\end{theorem}

\begin{demo}
We check that all invariant and ergodic probability measures $\mu$
with $\mu(\Sigma_M)=0$ are shadowed by $K(T,P)$, the kneading set.
This extends immediately to nonergodic measures, using the ergodic
decomposition and yields the result as $K(T,P)$ is finite (so that
$h_\top(\sigma,K(T,P))=0$. We shall finally prove the converse.

\step1{Almost every $A$ Markov at $0$ belongs to $\Sigma_M$}

Let $X_p$ is the set of $A\in\Sigma$ which are Markov at time $p$.
We claim that $X_p\subset X_{p+1}$.  Indeed, for $A\in\Sigma$,
$[\dots]_+$ denoting the symbolic cylinders in $\Sigma_+$ (i.e.,
the set of sequences in $\Sigma_+$ starting with a given word),
 \begin{equation}\label{eq-inclusion-Markov}
    \sigma^{n+1}[A_{p-n-1}\dots A_0]_+
     = \sigma^{n}[A_{p-n}\dots A_0]_+ \iff
     \sigma(A_{p-n-1}) \supset [A_{p-n}\dots A_0]_+
 \end{equation}
as $\sigma|[A_{-n}\dots A_0]$ is one-to-one. But
 \begin{multline*}
   \sigma(A_{p-n-1}) \supset [A_{p-n}\dots A_p]_+ \\
    \implies \sigma(A_{p+1-(n+1)-1}) \supset [A_{(p+1)-(n+1)}\dots A_p]_+
    \supset [A_{(p+1)-(n+1)}\dots A_{p+1}]_+
 \end{multline*}
proving the claim. On the other hand, $\sigma^{-1}(X_p)=X_{p+1}$
so $X_p$ and $X_{p+1}$ both have the same measure for any
invariant probability measure.

Thus $X_p\subset X_{p+1}$ implies $X_p=X_{p+1}$ modulo a set
negligible wrt all invariant probability measures so that:
 $$
     X_0 = \bigcap_{p\in\ZZ} X_p = \Sigma_M
 $$
modulo a negligible set wrt any invariant probability measure as
claimed.

 \step2{Geometric consequence of non-Markovianness}

By eq. (\ref{eq-inclusion-Markov}), $A$ is non-Markov at time $0$
implies on the interval that, for arbitrarily large integers $n$:
 $$
   T\left<A_{-n-1}\right> \not\supset \left<A_{-n}\dots A_0\right>
 $$
The right hand side is an interval so it is {\bf connected}
implying that
 $$
   (\partial T\left<A_{-n-1}\right>) \cap\left<A_{-n}\dots A_0\right>
    \ne \emptyset
 $$
This implies that $A_{-n-1}\dots A_0$ is the beginning of the
itinerary of one endpoint of $\left<A_{-n-1}\right>$, i.e., a
sequence in $K(T,P)$. Thus $A$ non-Markov at time $0$ implies that
$A$ is shadowed at time $0$ by $K(T,P)$.

 \step3{Shadowed sequences are non-Markov}

We prove the converse implication (up to finitely many periodic
sequences).

Let $A\in\Sigma(T,P)$ be a sequence which is shadowed at time $0$:
there are infinitely many integers $n\to\infty$ such that (*)
$A_{-n}\dots A_0$ is the beginning of $K\in K(T,P)$ (we may assume
$K$ to be fixed as $K(T,P)$ is finite). Assume also that $A$ is
Markov at time $0$: on the interval we have: $TA_{-n-1}\supset
\left<A_{-n}\dots A_0\right>$ for all large $n$ (note that this
holds whether or not $P$ is generating).

$K$ is the itinerary of, say, the left endpoint of $(c,c')\in P$.
Assume for simplicity that $T|(c,c')$ is increasing.

Take a large integer $n$ with both properties. Let $(d,e):=
\left<A_{-n}\dots A_0\right>$. By the the assumptions,
$TA_{-n-1}=T(c,c')\supset (d,e)$ and $Tc\in[d,e]$, so that $c$ is
the left endpoint of $\left<A_{-n-1}\dots A_0\right>$.

Consider $n'>n$ also satisfying the first property. Let $p=n'-n$.
$T^{p}$ maps $\left<A_{-n'}\dots A_0\right>$ to $\left<A_{-n}\dots
A_0\right>$ as $A$ is Markov at time $0$. Hence the endpoint $c$,
which is the left endpoint of both previous intervals, is
$p$-periodic. More precisely, for any $\eps>0$, there exists
$\delta>0$ such that $0<t<\delta\implies f^p(c+t)\in(c,c+\eps)$.
This implies that $(A_k)_{-n'\leq k\leq 0}$ is $p$-periodic.

Taking $A$ non-Markov for all times and $n'\to\infty$, we get that
$A$ is a periodic sequence from $K(T,P)$, finishing the proof of
the theorem.
\end{demo}

\begin{remark}
The non-Markov part of $\Sigma(T,P)$ is non-trivial for many
piecewise monotone maps. For instance, for an infinitely
renormalizable map, it carries the corresponding odometer. Even
for $\beta$-transformation, it can carries non-periodic measures,
for instance ones isomorphic to rotations (see
\cite{BuzziJohnson}). Note that this does not says that one cannot
find a Borel isomorphism between $\Sigma(T,P)$ and $\HS$, just
that it is not given by the natural map $\hat\pi$.
\end{remark}

\begin{exercise}
Find a piecewise monotone map $(T,P)$ with an invariant
probability measures $\mu$ on the \emph{one-sided} Markov shift
defined by $\MD$ which is not sent by $\hat\pi$ to an isomorphic
one on $\Sigma_+(T,P)$. \emph{Hint:} Use for instance that, for
piecewise invertible maps $T$ with invariant probability measure
$\mu$, the measure-theoretic Jacobian $d\mu\circ T/d\mu$ is an
invariant of measured conjugacy.
\end{exercise}

\begin{remark}\label{rem-exp-shad}
G. Keller has developped another approach to Hofbauer's
construction by trying to keep the smoothness of the interval map
in the picture. In particular he showed that, for "nice" interval
maps, like unimodal maps with negative Schwartzian derivative and
no attracting periodic orbits, the measures that cannot be lifted
to the Markov shift have been characterized by G. Keller and H.
Bruin \cite{KellerLifting,KellerBruin} as those having a zero or
negative Lyapunov exponent. Assumptions on Lyapunov exponents are
more general (recall Ruelle's inequality) but are also seem much
harder to generalize.
\end{remark}

\subsection{Entropy Bound}

\begin{proposition}\label{prop-shadow-entropy}
If $\mu$ is an invariant probability measure on $\Sigma$ shadowed
by $S\subset\Sigma_+$, then
 $$
   h(\sigma,\mu) \leq h_\top(S,\sigma)
 $$
\end{proposition}

Recall that
$h_\top(S,\sigma)=\limsup_{n\to\infty}\frac1n\log\#\{A_0\dots A_n:
A\in S\}$.

\begin{demo}
Observe that the canonical partition of $P^\ZZ$ (which we'll also
denote by $P$) is generating so that, according to Corollary
\ref{coro-SMB}, it is enough to find, for some $\lambda>0$, for
arbitrarily large $n\geq 1$, a small collection of $P,n$-cylinders
with a union of measure $\geq\lambda$.

Recall that $S\subset\Sigma_+$ and let $S_n:=\{A_0\dots
A_{n-1}:A\in S\}$. Let $\eps>0$. By definition of the topological
entropy, there exists $N_0<\infty$ such that for all $n\geq N_0$,
 $$
    \#S_n \leq e^{n(h_\top(\sigma,S)+\eps)}.
 $$
Increase $N_0$ if necessary so that the binominal coefficients
satisfy:
 $$
    \binom{2n/N_0}{n} \leq e^{\eps n}
 $$

 According to the main assumption the following measurable
function over $\Sigma$ is finite a.e.:
 $$
   n(A) := \min\{n\geq N_0: A_{-n+1}\dots A_0\in S_n\}
 $$
Thus there exists $M_0<\infty$ such that
 $$
   \mu(\{x\in \Sigma:n(x)>M_0\}) < \eps/\log\# P
 $$
By the ergodic theorem, there exists $N_1<\infty$ and $X_1$ such
that $\mu(X_1)>\lambda$ and for all $n\geq N_1$, for all $x\in
X_1$:
 $$
   \#\{0\leq k<n: n(\sigma^kx)\leq M_0 \} < \frac{\eps}{\log\#
   P} n
 $$

For $n\geq \max\left((\log\#P/\eps)M_0+N_1\right)$ we cover $X_1$
with a small number of $P,n$-cylinders.

For $x\in X_1$, we are going to define $c(x)=(a,b)$ where $a,b$
are integer sequences satisfying:
 $$
    0<a_r<b_r<a_{r-1}<b_{r_1}<\dots <a_1<b_1\leq a_0=n
 $$
and:
 \begin{enumerate}[(i)]
  \item for $1\leq s\leq r$, $\ell_s:=b_s-a_s+1\geq N_0$;
  \item $A_{a_s}\dots A_{b_s}\in S_{\ell_s}$;
  \item $\sum_{s=1}^{r} |b_s-a_s+1| >
  \left(1-2\frac{\eps}{\log\#P}\right)n$
 \end{enumerate}

We claim that this will prove:
 $$
   r(P,n,X_1,\lambda) \leq \binom{2n/N_0}{n} (\#P)^{n(\eps/\log\#P)+M_0}
   \exp n(h_\top(S)+\eps) \leq
   \exp n(h_\top(S)+4\eps)
 $$
and therefore the entropy bound as $\eps>0$ is arbitrary. We
conclude the proof of the proposition by proving this claim.
Define $c(x)=(a,b)$ over $X_1$ by setting: $a_0=n$ and then,
inductively:
 \begin{itemize}
   \item $b_{s+1}:=\max\{b< a_s:n(\sigma^bx)\leq M_0\}$;
   \item $a_{s+1}:=b_{s+1}-n(\sigma^{b_{s+1}}x)$
 \end{itemize}
and $r:=\max\{s\geq0:a_s\geq 0\}$. Observe that $c$ takes at most
$\binom{2n/N_0}{n}$ distinct values.

$x|\bigcup_{s=1}^r[a_s,b_s]$ can be described by specifying for
each $1\leq s\leq r$, $x|[a_s,b_s]\in S_{\ell_s}$: there are at
most $\exp n(h_\top(\sigma,S)+\eps)$ choices.

Let $I_*:=[0,a_r[\cup \bigcup_{s=1}^r ]b_s,a_{s-1}[$. For $k\in
I_*$, either $k\in]a_{r+1},b_{r+1}[$ (which has length at most
$M_0$) or $n(T^kx)> M_0$, therefore this subset of $[0,n[$ has
cardinality at most $M_0+(\eps/\log\# P)n$ by the choice of $M_0$.
Thus there are at most $(\#P)^{\# I_*}\leq (\#P)^{M_0+(\eps/\log\#
P)n}$ distinct choices for $x|I_*$.

Multiplying these numbers of choices we get the claim.
\end{demo}

\subsection{Weak Rank One}

The non-Markov part of the symbolic dynamics of a piecewise
monotone map is shadowed by a finite set. As the entropy of a
finite set is zero, Proposition \ref{prop-shadow-entropy} implies
that the entropy of any invariant probability measure carried by
this part is zero. Can one go further?

A motivation would be to find other type of arguments to control
$\Sigma\setminus\Sigma_M$, especially in higher dimensions where
the kneading set is replaced by a set with positive entropy.

Here, finiteness is much stronger than having zero entropy so it
is natural to ask whether one can say more about these measures:

\medbreak
 \centerline{\emph{What are the measures shadowed by
finitely many sequences?}} \medbreak

The good news is that this property can be formulated as a natural
generalization of the classical, ergodic theory notion of rank one
\cite{King}: an automorphism $T$  of a probability space $(X,P)$
has
\new{rank one} if, for every finite measurable partition $P$, for
every $\eps>0$ and $N<\infty$, there exists a finite word $\omega$
with length at least $N$ such that; for almost every $x\in X$, for
all $n\geq n_0(x,\eps,N)$, one can write the $n,P$-itinerary of
$x$ as $G^0\omega^1G^1\dots G^{r-1}\omega^r G^r$ where:
 \begin{itemize}
  \item the sum of the lengths of the gaps $G^s$ is at most $\eps n$;
  \item $|\omega^s|=|\omega|$ and $\sum_{s=1}^r \#\{0\leq
  i<|\omega|: \omega^s_i\ne\omega_i\}\leq\eps n$.
 \end{itemize}
Indeed, the propery of weak rank one is obtained by loosening the
"approximate copies of a long word" to "approximate copies of long
prefixes of a word":

\begin{definition}
An automorphism $T$  of a probability space $(X,P)$ has
\new{weak rank one} \cite{BuzziWRO} if, for every finite measurable partition $P$,
for every $\eps>0$ and $N<\infty$, there exists a finite word
$\omega$ such that; for almost every $x\in X$, for all $n\geq
n_0(x,\eps,N)$, one can write the $n,P$-itinerary of $x$ as
$G^0\omega^1G^1\dots G^{r-1}\omega^r G^r$ where:
 \begin{itemize}
  \item the sum of the lengths of the gaps $G^s$ is at most $\eps n$;
  \item $N\leq |\omega^s|\leq|\omega|$ and $\sum_{s=1}^r \#\{0\leq
  i<|\omega_s|: \omega^s_i\ne\omega_i\}\leq\eps n$.
 \end{itemize}
\end{definition}

\begin{proposition}
An automorphism $T$ of a probability space $(X,\mu)$ is weak rank
one iff, there exists a finite, measurable, generating partition
$P$ with respect to which $\mu$ is shadowed by a single sequence
\cite{BuzziWRO}.
\end{proposition}

Weak rank one is weaker than the usual generalizations of rank
one:

\begin{proposition}
All systems which have locally finite rank \cite{Ferenczi} are
weak rank one and the inclusion is strict. However there exists
zero entropy loosely Bernoulli systems which are not weak rank
one. In the $\bar d$-metric (see \cite{Rudolph, Shields}), the
generic system is \emph{not} weak rank one (i.e., the weak rank
one systems are contained in a countable union of closed subset
with empty interiors).
\end{proposition}

\begin{exercise}
Recall that system is of rank $r<\infty$ if it satisfies the
definition of rank one with the modification that instead of a
single word $\omega$ one has $r$, possibly distinct, words
$\omega(1),\dots,\omega(r)$. Show that such systems are weak rank
one.
\end{exercise}

The bad news is:

\begin{question}
Find a \emph{bona fide} ergodic property holding for all systems
with weak rank one, but not all some zero entropy systems.
\end{question}

\subsection{Shadowing Sequences}
A \new{shadowing sequence} is a sequence which shadows an
invariant probability measure. Not every sequence is shadowing:

\begin{exercise}
Show that any word appearing in a shadowing sequence must appear
infinitely many times.
\end{exercise}

\begin{question}
What can one say about the repetition times of $x$ (i.e.,
$n(\ell)\geq1$ minimum such that $x_{n(\ell)}\dots
x_{n(\ell)+\ell-1}=x_0\dots x_{\ell-1}$)? And about the return
times of $\mu$? Is there a relationship?
\end{question}

In fact, shadowing sequences are exceptional from a
measure-theoretic point of view:

\begin{exercise}
Show that the set of sequences in $\{0,1\}^\NN$ which shadow at
least one measure has zero measure wrt the
$(\frac12,\frac12)$-Bernoulli measure. \emph{Hint:} Use the
ergodic theorem to restrict the measures that could be shadowed.
\end{exercise}

Does the same occur in "concrete" situations:

\begin{question}
Is the set of $\beta>1$ such that the itinerary of $1^-$ is a
shadowing sequence a set of zero Lebesgue measure? (It is known
\cite{BuzziJohnson} that this set is uncountable).
\end{question}

One the other hand, it is common from the topological point of
view:

\begin{exercise}
Show that the set of sequences that shadow some measure is of
second Baire category (contains a dense countable intersection of
open subsets). \emph{Hint:} Consider, for each $n\geq1$,
$G_n:=\{x\in\{0,1\}^\NN:x=x_0\dots x_{n+m-1}(x_0\dots
x_{n-1})^p\dots$ for some $p\geq (n+m)^2\}$. See \cite{BuzziWRO}.
\end{exercise}

Some shadowing sequences $x$ satisfy:
$\mu_x:=\lim_{n\to\infty}\frac1n\sum_{k=0}^{n-1}
\delta_{\sigma^kx}$ exists in the weak star topology.

\begin{exercise}
Show that when the limit $\mu_x$ exists, $x$ can shadow only that
measure. Build two shadowing sequences one such that the measure
exists and one such that it does not.
\end{exercise}

Say that two measures are \emph{co-shadowed} if there exists a
sequence which shadows both of them.

\begin{exercise}
Show that two measures that are shadowed by two possible distinct
sequences are co-shadowed iff they have the same support.
Generalize this to a countable collection of measures
\cite{BuzziWRO}.
\end{exercise}

\begin{question}
Does there exists a sequence which shadows uncountable many
distinct measures? Does there exists a "universal shadowing
sequence", i.e., such that any shadowable measure with full
support is shadowed by it? Do these sequences form a set of second
Baire category?
\end{question}

\section{SPR Markov shifts}\label{sec-SPR}

The heart of the entropy theory of interval maps presented here is
the representation of their dynamics by combinatorial systems.
These are countable state Markov shifts but of a special kind,
called \emph{Strongly (or Stably) Positive Recurrent} Markov
shifts. They are the closest to the finite state Markov shifts:
most of the classical results generalize to them as we are going
to explain.

For countable Markov shift we refer to
\cite{Kitchens,Sarig0,Sarig1,RuetteClass} and especially to the
treatise \cite{GurevicSavchenko} which considers the SPR ones in
more details.

Recall the following definitions. A (countable state) \new{Markov
shift} is $\Sigma=\Sigma(G)$ where $G$ is a countable oriented
graph and $\Sigma:=\{x\in G^\ZZ:\forall p\in\ZZ\; x_p\to x_{p+1}$
on $G\}$ endowed with the left-shift $\sigma$. Note that
$\Sigma(G)$ is the vertex shift of $G$ - see \cite{LM} for the
edge shift). One can define an obvious one-sided Markov shift
$\Sigma_+(G)$ where the left-shift is again denoted by $\sigma$.

$G$ or $\Sigma$ are called
\new{irreducible} if there is a path from any vertex to any other
vertex in $G$ (i.e., $G$ is \emph{strongly connected}) (or,
equivalently, if $\Sigma$ is topologically transitive).

The \new{spectral decomposition} of a Markov shift is the at most
countable collection of irreducible Markov shifts defined by the
maximum subgraphs  of its defining graph which are strongly
connected components. Any ergodic invariant probability measure is
carried on a single piece of the spectral decomposition.

The \new{period} of an irreducible $G$ is the highest common
factor of all the lengths of its loops.

The \new{outdegree} of a vertex $v\in G$ is $\#\{w\in G:v\to w\}$.
The outdegree of $G$ is the supremum of the outdegrees of all
vertices.

\subsection{Definition}
Let $\Sigma$ be an irreducible Markov shift defined by a countable
oriented graph $G$. The basic classification of such objects from
our point of view is due to Vere-Jones \cite{VJ0,VJ1,VJ2}, which
remarked that one could extend the probabilistic classification of
stochastic matrices (related to the underlying Markov chain) to a
large class of positive matrices.

\begin{theorem}[Vere-Jones]\label{theo-VJ}
Fix an arbitrary vertex $v$ of $G$. Let $f_n$ be the number of
first returns after time $n$ to $v$. Let $\ell_n$ be the number of
returns after time $n$ to $v$. Let $R$ be the radius of
convergence of $\ell(z):=\sum_{n\geq1} \ell_nz^n$:
 $$
   R = \exp -h(\Sigma)
 $$
Assume that $R$ is finite. Then $\sum_{n\geq1} f_nR^{-n}\leq 1$.
\begin{enumerate}[(1)]
 \item $\sum_{n\geq1} f_nR^{-n}<1$: $G$ is \emph{transient} and $A$ admits no left or right
 positive eigenvectors;
 \item $\sum_{n\geq1} f_nR^{-n}=1$ and $\sum_{n\geq1}
 nf_nR^{-n}=\infty$: $G$ is \emph{null recurrent} and $A$ has both a left and a
 right eigenvectors, say $\ell$ and $r$, with $\sum_i
 \ell_i\cdot v_i=\infty$;
 \item $\sum_{n\geq1} f_nR^{-n}=1$ and $\sum_{n\geq1}
 nf_nR^{-n}<\infty$: $G$ is \emph{positive recurrent} and $A$ has both a left and a
 right eigenvectors, say $\ell$ and $r$, with $\sum_i
 \ell_i\cdot v_i<\infty$;
\end{enumerate}
Moreover, one can distinguish among the last type, those graphs
which are \emph{SPR} defined by the condition:
 $$
    \limsup_{n\to\infty} |f_n|^{1/n} < R.
 $$
In the above cases, the iterates of $A$ have the following
behavior as $n\to\infty$ for any (and then all) $i,j$:
 \begin{enumerate}[(1)]
  \item $\sum_{n\geq0} A^n_{ij}R^{-n}<\infty$;
  \item $\sum_{n\geq0} A^n_{ij}R^{-n}=\infty$ and $A^n_{ij}R^{-n}\to 0$;
  \item $A^n_{ij}R^{-n}\to d\cdot \ell_ir_j \ne 0$ where $d$ is the period of $G$ and $r$ and $\ell$ are positive eigenvectors: $\sum_i \ell_iA_{ij}=R\ell_j$ and
  $\sum_j A_{ij}r_j=R r_i$ normalized by $\sum_{k}
  \ell_k\cdot f_k=1$.
\end{enumerate}
In the SPR case, we additionally have:
 $$
   \limsup_{n\to_infty} \frac1n\log\left(A^n_{ij}R^{-n}-d\cdot \ell_i r_j\right) < 0
 $$
\end{theorem}

B.M. Gurevi\v{c} \cite{Gurevic} discovered the link with dynamics
and that there was in fact a probability measure hidden in these
positive matrices.

\begin{theorem}[Gurevi\v{c}]\label{theo-Gurevic}
In the Vere-Jones theorem, $R=\exp -h(G)$, where $h(G)$ is the
supremum of the measured entropies of $\Sigma$. This supremum is
achieved by at most one measure (called the \emph{maximum
measure}).

$G$ is positive recurrent if and only if there exists a maximum
measure. In this case, this measure $\mu_M$ is Markov and
satisfies:
 $$
   \mu_M([v_1\dots v_n]) = R^{-n} \ell_{v_1}r_{v_n} \text{ or $0$
   if } [v_1\dots v_n]=\emptyset.
 $$
\end{theorem}

\begin{exercise}
Show that the above data defines a probability measure $\mu_M$
which is invariant.
\end{exercise}

\begin{proposition}
Such measures as $\mu_M$ are finite extension of a Bernoulli
measure \cite{MaxBern}. The period of the extension is called the
period of the measure.
\end{proposition}

\begin{exercise}
Show that a positive recurrent graph is SPR iff the maximum
measure is exponentially filling: for any non-empty open
$U\subset\Sigma$:
 $$
   \limsup_{n\to\infty} \frac1n\log\mu\left(\bigcup_{0\leq k< n} \sigma^{-k}U\right)
     < 0
 $$
\end{exercise}

The SPR graphs have other characterizations:

\begin{proposition}[Gurevi\v{c}]
Let $G$ be irreducible with finite entropy. $G$ is SPR iff any
subgraph $G'$ obtained by removing one arrow satisfies:
$h(G')<h(G)$.
\end{proposition}

\begin{exercise}
Show that if $G$ is irreducible with finite entropy, then $G$ is
SPR iff any subgraph $G'$ obtained by removing one arrow is still
positive recurrent.
\end{exercise}

From our point of view the following characterization is
fundamental, showing that graphs that are "simple at infinity" are
good:

\begin{theorem}[Gurevi\v{c}-Zargaryan]
Let $G$ be an irreducible Markov shift with finite entropy. Define
its \new{entropy at infinity} by:
 $$
   h_\infty(G) := \inf_{F\subset\subset G\; \eps>0}
     \sup \{h(\sigma,\mu) : \mu\in\Prob(\Sigma) \text{ s.t. }
     \mu(\{x\in\Sigma:x_0\in F\})<\eps\}
 $$
where $F$ ranges over the finite subsets of vertices of $G$. $G$
is SPR iff $h_\infty(G)<h(G)$.
\end{theorem}

\subsection{Artin-Masur Zeta Function}

\begin{definition}
If $F$ is a subset of the vertices of the graph $G$, the
Artin-Masur \new{zeta function at $F$} is the formal power series:
 $$
   \zeta^G_F(z) := \exp \sum_{n\geq1} \frac{z^n}n
      \#\{ x\in\Sigma(G): \sigma^n x=x \text{ and }
      \{x_0,x_1,\dots\}\cap F\ne\emptyset\}
 $$
(provided each cardinality above is finite).

For finite $F\subset G$, $\zeta^G_F(z)$ is a \new{semi-local} zeta
function. For $F$ reduced to a single vertex $a$,
$\zeta^G_a(z):=\zeta^G_{\{a\}}(z)$ is called the \new{local} zeta
function.
\end{definition}

The full zeta function $\zeta^G_G(z)$ may fail to be defined or to
have the expected properties (see \cite[Thm. 9.4 and
following]{GurevicSavchenko}).

\begin{theorem}[Gurevich-Savchenko, Buzzi]\label{theo-MS-zeta}
Let $G$ be an irreducible countable oriented graph with
$h(G)<\infty$. Let $d$ be its period, i.e., the largest common
divisor of the length of all loops in $G$.

\begin{enumerate}[(1)]
 \item Each semi-local zeta function of $G$ defines a holomorphic
function over $|z|<\exp -h(G)$.
 \item These functions have a meromorphic extension to $|z|<\exp
-h_\infty(G)$.\cite{BuzziPQFT}
 \item $G$ is SPR iff any of its local zeta functions $\zeta^G_a$ has a meromorphic
extension to $|z|<R_a$ for $R_a>e^{-h(G)}$ without zeroes and
whose only singularities are simple poles at $z=e^{-h(G)} e^{2i\pi
k/d}$, for $k=0,\dots,d-1$.\cite{GurevicSavchenko}
\end{enumerate}
\end{theorem}

\subsection{Classification}

The following result shows that such systems are in some sense
very simple (or completely chaotic depending on your point of
view).

\begin{theorem}[Boyle-Buzzi-Gomez]
To each irreducible Markov shift $\Sigma$ associate its entropy
$h(G)$ and the list of the periods of its maximum measures
(repeated according to multiplicities).

This data is a complete invariant for entropy-conjugacy among SPR
Markov shifts.
\end{theorem}

\begin{remark}
In fact a stronger notion of conjugacy is established in
\cite{BBG}: one finds in each graph a vertex and a Borel conjugacy
between the sets of all sequences that visit infinitely often in
the past and in the future these vertices, except for those that
remain in two finite subgraphs and a (necessarily countable) set
of periodic sequences.
\end{remark}

\begin{remark}
For subshifts of finite type this is a classical result of Adler
and Marcus \cite{AdMa}.
\end{remark}

The idea of the proof of the above is quite simple. Let us sketch
it. Consider two SPR graphs. One can assume these graphs to be
loop graphs with many first return loops of each large length
after discarding some subshifts of finite type (which have
strictly smaller entropy by the SPR property). Then one modifies
these graphs by removing loops of increasing length from one or
the other until they both have the same number of loops of all
lengths. These removals are operated on first return loops. A
first return loop $\ell$ is removed and first return loops
$L\ell^n$, $n\geq1$, for any first return loop $L\ne\ell$ are
added back so that no other sequence has been removed from the
shift. The Markov shift obtained in this limit is Borel
isomorphic, except for the set of removed periodic sequences, to
the initial ones \emph{if}, for instance,  one has been able to
remove only first return loops that were already present in these
initial graphs. The SPR assumption ensures that the number of
loops to be removed is so small that the previous condition can
easily be fulfilled (see Theorem \ref{theo-MS-zeta}).

\section{Application to Piecewise Monotone
Maps}\label{sec-application}

Let $(T,P)$ be a piecewise monotone map with its symbolic dynamics
$\Sigma_+(T,P)$ and the natural extension of it, $\Sigma(T,P)$.

By Theorem \ref{theo-smallentropy1}, $\Sigma_(T,P)$, after
discarding a subset carrying only zero entropy measures and
finitely many periodic sequences, is Borel isomorphic to the
Markov shift $\HS$ defined by its Markov diagram $\MD$.

We claim that $h_\infty(\MD)=0$ so that, for every $H>0$, the
spectral decomposition of $\HS$ contains only finitely many
irreducible Markov shifts with entropy $\geq H$ and that each of
these is SPR.

This follows from the description of Theorem \ref{theo-PMM-MD}.
$\MD$ has finite outdegree (bounded by $\#P$) and, after removing
the finite subset
 $$
   \MD_N:=\{\min(K_0\dots K_n): K\in K(T,P) \text{ and }n\leq N\}
 $$
the remaining graph $\MD\setminus\MD_N$ has outdegree at most $2$
and the following property: if $\alpha_0\to\dots\alpha_n$ is a
path on this graph with $\alpha_0$ and $\alpha_n$ with outdegree
$>1$ in $\MD\setminus\MD_N$, then $n\geq N$.

\begin{exercise}
Check that the above property holds and that it implies that
$h_\infty(\MD)=0$.
\end{exercise}

\subsection{Measures of Maximum or Large Entropy}

We prove in this section the following:

\begin{theorem}[Hofbauer, Boyle-Buzzi-Gomez]
A piecewise monotone map with nonzero entropy has a finite number
of ergodic measure, each of which is up to a period, a Bernoulli
probability measure \cite{HofbauerPMM,MaxBern}.

More precisely, the topological entropy together with the list
with repetitions of the periods of the maximum measure is a
complete invariant for the entropy-conjugacy of the natural
extensions of these maps \cite{BBG}.
\end{theorem}

\begin{remark}
This is false if one removes the words "natural extensions of".
\end{remark}

\begin{question}
Is a piecewise monotone map with nonzero topological entropy
always entropy-conjugate to \emph{some} countable state Markov
shift?
\end{question}

The first part of the theorem follows from Gurevi\v{c} result that
a Markov shift has exactly one maximum measure (ergodic, invariant
probability measure with maximum entropy) on each piece of the
spectral decomposition which is positive recurrent. Moreover these
measures are Markov and this implies by \cite{MaxBern} that they
are Bernoulli up to a period.

\begin{exercise}
Show that there is a bijection which is entropy and ergodicity
preserving between a dynamical system and its natural extension.
\end{exercise}

Hence the above implies that piecewise monotone maps with nonzero
topological entropy have a finite number of maximum measures.

By exploiting the very special structure of the Markov diagram of
a piecewise monotone map one can obtain quantitative bounds.

\begin{theorem}[Hofbauer,Buzzi]
A unimodal map with positive topological entropy has exactly one
maximum measure \cite{HofbauerUNI}.

A piecewise monotone map $T:[0,1]\to[0,1]$ with nonzero entropy
and $N$ pieces under the normalization $T(\{0,1\})\subset\{0,1\}$
has at most $4(N-1)$ maximum measures (only $3(N-1)$ if it is
continuous) \cite{BuzziNES}.
\end{theorem}

However the following very natural conjecture remains open:

\begin{question}
A piecewise monotone map with nonzero entropy and $N$ pieces has
at most $N-1$ maximum measures.
\end{question}

\begin{remark}
Another approach to this problem is to use the following results:
(1) for a piecewise monotone map with slope constant in absolute
value, the maximum measures are exactly the ergodic absolutely
continuous invariant probability measures (this is not totally
obvious because of the discontinuity \cite{DKU}); (2) any
piecewise monotone map with nonzero entropy admits as a
topological factor a piecewise monotone map of the previous type
\cite{MT}; (3) the number of distinct ergodic absolutely
continuous invariant probability measure of, say, an  interval map
with piecewise constant slope and $N$ pieces is bounded by $N-1$.
The problem is to understand what measures are collapsed under the
factor map and to relate the number of maximum measures to the
successive renormalizations of the map.
\end{remark}

\subsection{Periodic points}

The above estimate $h_\infty(\MD)=0$ yields, using Theorem
\ref{theo-MS-zeta}:

\begin{theorem}[Milnor-Thurston, Hofbauer]
Let $(T,P)$ be a piecewise monotone map with non zero entropy. The
Artin-Masur zeta function:
 $$
   \zeta_T(z):=\exp\sum_{n\geq1} \frac{z^n}n
   \#\{A\in\Sigma_+(T,P):\sigma^nA=A\}
 $$
is holomorphic on $|z|<e^{-h_\top(T)}$. It has a meromorphic
extension to $|z|<1$ \cite{MT} (see \cite{Ruelle2} for a later,
simpler proof).

The periodic orbits are equidistributed according to the maximum
measures:
 $$
   \frac1{\#\{A\in\Sigma_+(T,P):T^nA=A\}} \sum_{A\in\Sigma_+(T,P)\text{ s.t. }A=T^nA}
     \delta_A \longrightarrow^{n\to\infty}_{n\in p\NN} \frac1{\sum_{i=1}^rp_i}\sum_{i=1}^r p_i\mu_i
 $$
where $p:=\operatorname{lcm}(p_1,\dots,p_r)$, $\mu_1,\dots,\mu_r$
are the maximum measures on $\Sigma_+(T,P)$ and $p_i$ is the
period of $\mu_i$.
\end{theorem}

\section{Further Results}\label{sec-conclusion}

\subsection{Generalization to Other Dynamics}

The results presented above are well-known in the classical
hyperbolic setting, e.g., for subshifts of finite type \cite{LM}.
They in turn can be vastly generalized to various non-uniformly
expanding systems.

To begin with, one can define a class of subshifts over finite
alphabets, the
\new{subshifts of quasi-finite type} \cite{BuzziQFT}, which includes
both the subshifts of finite type and the symbolic dynamics of
piecewise monotone maps and satisfy most of the properties
explained in this paper:

\begin{definition}
A subshift $\Sigma$ over a finite alphabet is of \new{quasi-finite
type} if the minimum words of length $n$, $\mathcal M_n$, satisfy:
 $$
    \limsup_{n\to\infty} \frac1n\log\#\mathcal M_n <
    h_\top(\sigma|\Sigma).
 $$
\end{definition}

This class also contains new dynamics, e.g., multidimensional
$\beta$-transformations, i.e., self-maps of $[0,1]^d$ defined by
$x\mapsto A.x \mod \ZZ^d$ where $A$ is an affine transformation
such that $\|A^N.v\|>\|v\|$ for some $N$ and all
$v\in\RR^d\setminus\{0\}$. The entropy theory of arbitrary
piecewise affine maps in higher dimension is mostly unknown
\cite{Kruglikov} beyond similar maps and piecewise affine surface
homeomorphisms \cite{BuzziPWAH}.

One can further extend the approach explained here to
\new{entropy-expanding} maps, which are multidimensional smooth
maps with critical points such that the entropy of smooth
submanifolds is smaller than their topological entropy. This class
includes arbitrary $C^\infty$ interval maps with nonzero entropy
\cite{BuzziSIM,BuzziRuette} and maps like:
 $$
    (x,y)\mapsto (1.9-x^2+\eps y,1.8-y^2+\eps x)
 $$
for small $|\eps|$. They can be treated either geometrically
\cite{BuzziSMF,BuzziEE} or, under a "general position assumption",
by symbolic techniques with the introduction of \new{puzzles of
quasi-finite type} \cite{BuzziPQFT} which generalize to
combinatorial puzzles (as defined in complex dynamics) the notions
above. One obtains in this way, e.g., that they possess finitely
many maximum measures with Bernoulli extensions (up to a period).

\subsection{Thermodynamical Formalism}

The \new{thermodynamical formalism} introduced in the uniformly
hyperbolic setting by Sinai and Ruelle involves the introduction
of a weight in the countings defining the entropy. This
generalization of entropy is called "pressure" (it is in fact the
free energy per site of a one-dimensional crystal formally
associated to the dynamical system) \cite{Ruelle1, Walters,
KellerEquilibrium}. The measures maximizing the pressure are
called the equilibrium measures.

This generalization is very important for several reasons. For
instance, if $f$ is a smooth uniformly expanding map of the
circle, the equilibrium measure for the weight $|f'|$ is the
absolutely continuous invariant measure. Most of the entropy
theory can generalized to this setting either when in
one-dimension (\cite{BaladiKeller,Ruelle2} and the references in
the book \cite{BaladiBook}) or when there is enough uniformity in
the expansion (see, e.g.,
\cite{BuzziThermo,BuzziKeller,BuzziConfo}) so that the regularity
properties of the weight (bounded variation or Holder continuity)
also hold in the symbolic representation.

\subsection{Further Questions}

We have not been able to discuss \emph{stability issues}, that is,
the (dis)continuity of the map that associates to a piecewise
monotone map its set of maximum measures (see
\cite{RaithStability}).

We can only refer the reader to the already well-developed
\emph{dimension theory} \cite{PesinBook,Hofbauer-D} which in
particular relates Hausdorff dimension, entropy and Lyapunov
exponents or the developping \emph{symbolic extension theory}
\cite{BFF,BD,ND} which has already obtained deep refinements of
the local entropy estimates presented in section
\ref{sec-h-interval}.

\end{document}